\pdfoutput=1
\RequirePackage{ifpdf}
\ifpdf 
\documentclass[pdftex]{sigma}
\else
\documentclass{sigma}
\fi

\usepackage{enumitem}

\newcommand{\R}[1]{\mathbb{R}^{#1}}
\newcommand{\dd}{{\rm d}}
\newcommand{\ii}{\mathbf{i}}

\numberwithin{equation}{section}

\newtheorem{Theorem}{Theorem}[section]
\newtheorem*{Theorem*}{Theorem}
\newtheorem{Corollary}[Theorem]{Corollary}
\newtheorem{Lemma}[Theorem]{Lemma}
\newtheorem{Proposition}[Theorem]{Proposition}
 { \theoremstyle{definition}
\newtheorem{Definition}[Theorem]{Definition}

\newtheorem{Example}[Theorem]{Example}
\newtheorem{Remark}[Theorem]{Remark} }

\begin{document}

\allowdisplaybreaks

\newcommand{\arXivNumber}{2103.00458}

\renewcommand{\PaperNumber}{038}

\FirstPageHeading

\ShortArticleName{Geometrical Aspects of the Hamiltonization Problem of Dynamical Systems}

\ArticleName{Geometrical Aspects of the Hamiltonization Problem\\ of Dynamical Systems}

\Author{Misael AVENDA\~NO-CAMACHO, Claudio C\'esar GARC\'IA-MENDOZA, \newline Jos\'e Crisp\'in RU\'IZ-PANTALE\'ON and Eduardo VELASCO-BARRERAS}
\AuthorNameForHeading{M.~Avenda\~no-Camacho et al.}

\Address{Departamento de Matem\'aticas, Universidad de Sonora, M\'exico}
\Email{\href{mailto:misael.avendano@unison.mx}{misael.avendano@unison.mx}, \href{mailto:a214200511@unison.mx}{a214200511@unison.mx}, \href{mailto:jose.ruiz@unison.mx}{jose.ruiz@unison.mx}, \newline
\hspace*{14mm}\href{mailto:eduardo.velasco@unison.mx}{eduardo.velasco@unison.mx}}

\ArticleDates{Received March 02, 2021, in final form May 10, 2022; Published online May 20, 2022}

\Abstract{Some positive answers to the problem of endowing a dynamical system with a~Hamiltonian formulation are presented within the class of Poisson structures in a geometric framework. We address this problem on orientable manifolds and by using decomposable Poisson structures. In the first case, the existence of a Hamiltonian formulation is ensured under the vanishing of some topological obstructions, improving a~result of Gao. In the second case, we apply a variant of the Hojman construction to solve the problem for vector fields admitting a transversally invariant metric and, in particular, for infinitesimal generators of proper actions. Finally, we also consider the hamiltonization problem for Lie group actions and give solutions in the particular case in which the acting Lie group is a low-dimensional torus.}


\Classification{37J06; 37J39; 53D17; 37C86; 70G45; 37C79}

\section{Introduction}

The hamiltonization problem on a smooth manifold is the question of whether a dynamical system admits a Hamiltonian formulation. In the framework of Poisson Geometry, given a~vector field $X$ on a smooth manifold $M$, the problem consists in finding a scalar function~$h$ and a Poisson structure~$\pi$ on~$M$ such that
 \begin{equation}\label{EcIntroHam}
 X = \pi(\dd h, \cdot ).
 \end{equation}
If such a Poisson structure $\pi$ and function $h$ exist, then $X$ is said to be hamiltonizable on~$M$.

The problem of finding a Poisson structure with respect to which a flow of ODEs, or a~vec\-tor field, is Hamiltonian has received considerable attention from several decades ago. At~the earlier stages, this problem concerned to the existence of a Hamiltonian formulation of a~set of ODEs and their explicit formulation, mainly in two and three dimensions \cite{Abar-87,CairoFeix-92,Gumral-Nut,Hern-97,Whittaker}. Some geometric approaches to the hamiltonization problem can be found in the works of Whittaker~\cite{Whittaker}, Perlick~\cite{Perlick-92}, S.~Hojman \cite{Hoj-91,Hojman-96} and R.~Alvarado Flores et al.~\cite{AHA-06}. A solution to this problem is given by the Lie--K\"{o}nigs theorem \cite{Whittaker}, which requires the solvability of the ODEs. V.~Perlick gives a generalization of the Lie--Konigs theorem in a more global setting, but only when the original system is even-dimensional, so that a time-dependent symplectic structure exists~\cite{Perlick-92}. S.~Hojman solves the hamiltonization problem of a vector field by constructing a Poisson structure, provided that it admits an infinitesimal symmetry and a first integral~\cite{Hojman-96}. A geometric reformulation of Hojman results is presented in~\cite{AHA-06}, where the hamiltonization problem is solved for a vector field being in a two-dimensional Lie subalgebra and also admitting a first integral. In recent works, the hamiltonization problem has also been considered, both theoretically~\cite{Kozlov-2020} and in applied contexts \cite{Ballesteros_2020}.

In this paper, we present a geometric approach to the hamiltonization problem for vector fields (admitting zeroes, in general). In particular, we give positive answers to the hamiltonization problem for vector fields on orientable manifolds with sufficient first integrals, and also for vector fields admitting first integrals that are not necessarily regular. In the latter case, this improves the Hojman construction and the approach of~\cite{AHA-06}. Furthermore, we provide solutions to the hamiltonization problem for Lie group actions in some particular cases.

Our first result deals with vector fields on an $m$-dimensional manifold admitting $m-1$ first integrals. We show that such vector fields are Hamiltonian with respect to $m-1$ Poisson structures defined on the open set where the first integrals are independent.
More precisely, in Theorem \ref{teo:IntegrableHam} we prove that if a vector field~$X$ has $m-1$ first integrals $h_1,h_2,\dots, h_{m-1}$ functionally independent on a open dense set $U$, then there exist $m-1$ Poisson structures $\pi_1,\dots, \pi_{m-1}$ on $U$ with rank at most two such that
 \begin{equation*}
 \pi_i^\sharp(\dd h_j)=\delta_{ij}X,\qquad
 i,j=1,2,\dots, m-1,
 \end{equation*}
where $\delta_{ij}$ denotes the kronecker delta. That is, for each $i=1,2,\dots, m-1$, $X$ is Hamiltonian with respect to $\pi_i$ with Hamiltonian functions $h_i$ and the rest of first integrals are Casimir fuctions of $\pi_i$. Furthermore, the Poisson structures $\pi_{i}$ commute with respect to the Schouten--Nijenhuis bracket: $[\pi_i,\pi_j]=0$ for all $i,j=1,\dots,m-1$. These facts generalize the results presented in~\cite{Gumral-Nut}. We also describe the dependence of the Poisson structures on the choice of the volume form, which allows us to show in Proposition \ref{prop:IntegrableHam} that one can drop the orientability hypothesis. Moreover, we also improve the previous result by relaxing the hypothesis on the number of first integrals: if a vector field $X$ on an $m$-dimensional manifold admits $m-2$ independent first integrals on an open dense set and an invariant volume form, then $X$ is hamiltonizable under a suitable topological condition on the common level sets of the first integrals. This is the content of Theorem \ref{teo:VolHam}, which in the three-dimensional case recovers the hamiltonization criteria given by Gao in~\cite{Gao-2000} for Lotka--Volterra systems; but the topological hypotheses of our theorem were obviated there. These hypotheses are necessary in general as we exhibit it in Examples~\ref{ex:IntSing} and~\ref{exm:37}. On the other hand, by considering a correspondence between multivector fields and differential forms on orientable manifolds, we have formulated a criterion of hamiltonization by unimodular Poisson structures in terms of the existence of integrating factors of primitives (Theorem~\ref{teo:UnimodHam}). In this case, the given vector field is a modular vector field. Finally, Theorem~\ref{teo:DoubleFol} is a more general version of Theorem~\ref{teo:VolHam} that provides a hamiltonization criteria by means of differential forms that are not necessarily product of exact 1-forms: if a vector field preserves a leaf-wise volume of an oriented regular foliation of dimension $r$, and belongs to the kernel of a nowhere vanishing leaf-wise closed $(r-2)$-form, then the hamiltonization problem reduces to the triviality of the foliated de Rham cohomology in degree one for the foliation integrating the kernel of the closed form.

We also present various results on the hamiltonization problem where orientability is no longer required. Instead, we aim for decomposable Poisson structures. More precisely, in Theorem~\ref{teo:NormalSection}, we state necessary and sufficient conditions under which a given vector field $X$ on a~manifold~$M$ is hamiltonizable by a Poisson structure of the form $\pi = Y \wedge X$, for some vector field~$Y$ on~$M$. In particular, we recover the results of \cite{AHA-06,Hojman-96}, where $X$ must admit a regular first integral.
Moreover, following~\cite{LPS-12}, we present in Theorem \ref{teo:DiverInt} a hamiltonization criteria for vector fields admitting an infinitesimal symmetry and an invariant volume form. On the other hand, in Theorem \ref{teo:RegFol} we show that a vector field is hamiltonizable if it admits a first integral and a suitable 2-dimensional foliation. This is a generalization of~\cite[Theorem~2]{FGS-05}, where a similar result is proven for the case of nowhere vanishing vector fields with periodic flow. We also present geometric settings in which Theorem~\ref{teo:NormalSection} is applicable: in Theorem~\ref{teo:TrInvM} we show that a vector field is hamiltonizable under suitable conditions involving transversally invariant Riemannian metrics, and in Theorems~\ref{teo:HorInv} and \ref{teo:VertHam} we give hamiltonization criteria for vector fields tangent to the fibers of a submersion. These approaches allowed us to derive Theorem~\ref{teo:Periodic}, where we have shown that the positive answer to the hamiltonization problem for nowhere vanishing vector fields with periodic flow is given by a~topological condition: the non-compacity of the manifold. Also, we have derived a hamiltonization criteria for vector fields inducing a proper $\mathbb{R}$-action: such a vector field is hamiltonizable if and only if its orbit space is not compact (Theorem~\ref{teo:ProperHam}). Moreover, in Propositions~\ref{prop:InfHam} and~\ref{prop:InfHam2} we have applied Theorem~\ref{teo:TrInvM} to obtain a hamiltonization criteria for infinitesimal generators of proper actions of general Lie groups.\looseness=-1

We note that the Poisson structures in Sections~\ref{sec:Hamiltonization problem on orientable manifolds} and~\ref{sec:HamGood} are of rank at most two. This class of Poisson structures have properties that are generally not true for those of higher rank, and that helped us to derive the results of these two sections: they are conformal invariant, are in correspondence to orientable 2-dimensional foliations and are multiple of their nowhere vanishing Hamiltonian vector fields. Moreover, the set of singular points of such structures admits a simple description.

Finally, we consider the \emph{hamiltonization problem of Lie group actions}: given a Lie group action on a manifold $M$, whether a Poisson structure on $M$ exists with respect to which the action is Hamiltonian. We want to comment that we have not found in the literature any formulation nor results concerning to the hamiltonization problem of Lie group actions. For the abelian Lie group $\mathbb{T}^k$, we give some conditions under which there exists a bivector field $\pi$ on $M$ such that the infinitesimal generators are of the form \eqref{EcIntroHam}. If $\pi$ is Poisson, then the $\mathbb{T}^k$-action is Hamiltonian. In Theorem~\ref{teo:2Torus} we provide sufficient conditions to solve the hamiltonization problem for 2-dimensional torus actions. The one-dimensional case is the content of Theorem~\ref{teo:Periodic2}.

\section{Preliminaries} 

A \emph{Poisson structure} on a smooth manifold $M$ is a bivector field $\pi \in \Gamma \big({\wedge}^2 TM\big)$ satisfying the Jacobi identity
 \begin{equation*}
 [\pi,\pi] = 0.
 \end{equation*}
Here, the bracket $[-,-]$ stands for the Schouten--Nijenhuis bracket for multivector fields \cite[Section~1.8]{DZ}. The pair $(M,\pi)$ is called a \emph{Poisson manifold}.

The \emph{rank} of a Poisson structure $\pi$ at a point $x \in M$ is defined by $\operatorname{rank}\pi_{x} := \dim\pi^{\sharp}(T^{\ast}_{x}M)$, where the vector bundle map $\pi^{\sharp}\colon T^{*}M \to TM$ is given by the usual contraction $\pi^{\sharp}\alpha := \ii_{\alpha}\pi$, for all $\alpha \in T^{\ast}M$. A \emph{singular point} $x$ of $\pi$ is characterized by the condition that the rank of~$\pi$ is not constant around~$x$. Otherwise, the point $x$ is said to be a \emph{regular point} of~$\pi$.

The image $\pi^{\sharp}(T^{\ast}M)$ is an integrable distribution on $M$, called the \emph{characteristic distribution} of $\pi$, and each integral submanifold carries a symplectic structure canonically induced by the Poisson structure $\pi$ on $M$. Consequently, the integral symplectic submanifolds define a smooth symplectic foliation of $M$, which may be singular in general.

A vector field $X$ on a Poisson manifold $(M,\pi)$ is \emph{tangent} to the symplectic foliation if $X_{x} \in \pi^{\sharp}(T^{*}_{x}M)$, for all $x \in M$. Note that for every $\alpha \in \Gamma(T^{*}M)$ the vector field $X = \pi^{\sharp}\alpha$ is tangent to the symplectic foliation. In particular, if $\alpha = \dd{h}$, for some $h \in C^{\infty}(M)$, then $X$ is called \emph{Hamiltonian vector field}. In this case, the function $h$ is said to be a \emph{Hamiltonian function} for~$X$. We denote by $\mathfrak{ham}(M,\pi)$ the Lie algebra of all Hamiltonian vector fields on $M$.

The \emph{Casimir functions} of $\pi$ are the Hamiltonian functions $c \in C^{\infty}(M)$ of the zero vector field, that is, $\pi^{\sharp}\dd{c} = 0$. These functions are constant along the leaves of the symplectic foliation of $M$. Thus, every tangent vector field is also tangent to the level sets of every Casimir function. In particular, this is true for Hamiltonian vector fields.

The infinitesimal Poisson automorphisms of $\pi$, or \emph{Poisson vector fields}, for short, are the vector fields $X \in \Gamma(TM)$ such that ${\rm L}_{X}\pi = 0$. By $\mathfrak{poiss}(M,\pi)$ we denote the Lie algebra of Poisson vector fields, and the Lie subalgebra of tangent Poisson vector fields by $\mathfrak{poiss}_{{\rm tan}}(M,\pi)$. We recall that the Lie algebra $\mathfrak{ham}(M,\pi)$ is an ideal of both $\mathfrak{poiss}(M,\pi)$ and $\mathfrak{poiss}_{{\rm tan}}(M,\pi)$.

Now, we observe that if $\pi$ is a Poisson structure and $f$ is an arbitrary function on $M$, then~$f\pi$ is not a Poisson structure, in general. We say that $\pi$ is \emph{conformally invariant} if $f\pi$ is again a Poisson structure on $M$ for any $f \in C^{\infty}(M)$. A simple computation shows that
 \begin{equation*}
 [f\pi, f\pi] = -2f\,\pi^{\sharp}\dd{f} \wedge \pi.
 \end{equation*}

\begin{Lemma}\label{lemma:conformally}
Every Poisson structure with rank at most two is conformally invariant.
\end{Lemma}

This lemma follows from the fact that the 3-vector field $\pi^{\sharp} \dd{f} \wedge \pi$ is a section of the bundle $\wedge^3(\pi^{\sharp}(T^{*}M))$, which is zero due to our rank hypothesis.

{\bf Tangential Poisson cohomology.} Every Poisson structure $\pi$ on $M$ induces a cochain complex, called the \emph{Lichnerowicz--Poisson complex}. Its cohomology, denoted by $H^{\bullet}(M,\pi)$, is called the \emph{Poisson cohomology} of the Poisson manifold $(M,\pi)$ \cite{Lichnerowicz}. It is well-known that $H^{0}(M,\pi)$ consists of the Casimir functions, the 1-coboundaries are the Hamiltonian vector fields and the 1-cocycles are the Poisson vector fields of $(M,\pi)$ \cite{DZ,Lichnerowicz,Va-94,VK-88,VK-90}. So, the cohomology in degree one is the Lie algebra
 \begin{equation*}
 H^{1}(M,\pi) = {\mathfrak{poiss}(M,\pi)} \big/ {\mathfrak{ham}(M,\pi)}.
 \end{equation*}
The \emph{tangential Poisson cohomology} in degree one is the Lie subalgebra $H^{1}_{{\rm tan}}(M,\pi) \subseteq H^{1}(M,\pi)$ consisting of the Poisson cohomology classes with tangent representatives. More precisely,
 \begin{equation*}
 H^{1}_{{\rm tan}}(M,\pi) := {\mathfrak{poiss}_{{\rm tan}}(M,\pi)} \big/ {\mathfrak{ham}(M,\pi)}
 \end{equation*}
(see \cite[Definition~2]{AG}). Now, recall that the \emph{foliated de Rham complex} of a regular foliation consists of the graded algebra of foliated differential forms endowed with the foliated exterior derivative. Its cohomology is the \emph{foliated de Rham cohomology}.
Then, we have the following fact \cite[Chapter~5]{Va-94}.

\begin{Proposition}\label{prop:PoissTan}
The tangential Poisson cohomology in degree one of a regular Poisson manifold is isomorphic to the foliated de Rham cohomology of its symplectic foliation and, in particular, independent of the leaf-wise symplectic form.
\end{Proposition}

In the case when the symplectic foliation is \emph{simple}, that is, given by the fibers of a submersion, we have the following criterion for the vanishing of the tangential Poisson cohomology in degree one.

\begin{Proposition}\label{prop:H1triv}
Let $(M,\pi)$ be a regular Poisson manifold such that its symplectic foliation is given by the fibers of a submersion. If the fibers are connected and simply connected, then $H^{1}_{{\rm tan}}(M,\pi)=0$.
\end{Proposition}

This is consequence of Proposition \ref{prop:PoissTan} and of the fact that the topology of the fibers imply the triviality of the foliated de Rham cohomology in degree one \cite[Proposition~7.4]{DazHec-91}.

{\bf Orientable Poisson manifolds.} Recall that on each oriented $m$-dimensional manifold $M$, equipped with a volume form $\Omega$, there exists a one-to-one correspondence between $(m-2)$-differential forms and bivector fields given by the following formula:
 \begin{equation}\label{EcPiRHo}
 \ii_{\pi}\Omega = \varrho, \qquad \pi \in \Gamma\big({\wedge}^{2}TM\big),\quad \varrho \in \Gamma\big({\wedge}^{m-2}T^{\ast}M\big).
 \end{equation}
Here, the interior product of multivector fields and differential forms is defined by the rule $\ii_{A \wedge B} = \ii_{A} \circ \ii_{B}$, for any $A,B \in \Gamma(\wedge^{\bullet}TM)$.

Consider the $(m-2)$-differential form $\varrho = \frac{1}{2} \sum_{i,j} \varrho_{ij}\dd{x^{1}} \wedge \cdots \wedge \widehat{\dd{x^{i}}} \wedge \cdots \wedge \widehat{\dd{x^{j}}} \wedge \cdots \wedge \dd{x^{m}}$, where the 1-forms with $\ \widehat{}\ $ are omitted. We define the \emph{rank} of $\varrho$ at a point $x \in M$ as the rank of the matrix $[\varrho_{ij}(x)]$. Here, the set $\big\{x^{1},\dots,x^{m}\big\}$ is a local coordinate system on~$M$. We also recall that an orientable Poisson manifold is said to be \emph{unimodular} if admits a volume form invariant under the all Hamiltonian flows~\cite{We-97}.

\begin{Lemma}\label{lemma:RhoRank2}
If $\varrho$ is a closed $(m-2)$-differential form of rank at most two, then the bivector field $\pi$ defined by the formula \eqref{EcPiRHo} is a unimodular Poisson structure on $M$ of rank at most two.
\end{Lemma}

\begin{Lemma}\label{lemma:HamRho}
Let $\varrho$ be a $(m-2)$-differential form and suppose that the bivector field $\pi$ defined by the formula \eqref{EcPiRHo} is a Poisson structure on $M$. Then, a vector field $X$ on $M$ is a Hamiltonian vector field, with respect to $\pi$, if and only if there exists $h \in C^{\infty}(M)$ such that
 \begin{equation*}
 \ii_{X}\Omega = \dd{h} \wedge \varrho.
 \end{equation*}
In particular, a function $K \in C^{\infty}(M)$ is Casimir for $\pi$ if and only if $\dd{K} \wedge \varrho = 0$.
\end{Lemma}

Finally, one can verify that the Poisson bivector field given by \eqref{EcPiRHo} is unimodular if and only if $\varrho$ admits an (non-zero) integrating factor $a \in C^{\infty}(M)$, $\dd(a\varrho)=0$.

\section{Hamiltonization problem on orientable manifolds}
 \label{sec:Hamiltonization problem on orientable manifolds}

We begin by observing that every vector field on an orientable smooth $m$-dimensional manifold admitting $m-1$ independent first integrals is hamiltonizable.

\begin{Theorem}\label{teo:IntegrableHam}
Let $X$ be a vector field on an orientable $m$-dimensional manifold $M$. Suppose that $X$ admits $m-1$ first integrals $h_1, \dots, h_{m-1} \in C^{\infty}(M)$ which are independent on a open dense subset $U \subseteq M$. Then, there exist unique Poisson structures $\pi_{1}, \dots, \pi_{m-1}$ on $U$ of rank at most two such that
 \begin{equation}\label{HamX}
 \pi_{i}^{\sharp}\dd h_{j} = \delta_{ij}X, \qquad \text{for all} \quad i,j = 1, \dots, m-1.
 \end{equation}
Moreover, we have $[\pi_{i},\pi_{j}]=0$, for all $i$ and $j$.
\end{Theorem}
\begin{proof}
Fix a volume form $\Omega$ on $M$ and set $\beta := \dd h_1 \wedge \cdots \wedge \dd h_{m-1}$. First, define bivector fields~$\psi_{i}$ on $M$ by the formula
 \begin{equation}\label{eq:Psi-i}
 \ii_{\psi_{i}}\Omega = \dd h_{1} \wedge \cdots \wedge \widehat{\dd h_{i}} \wedge \cdots \wedge \dd h_{m-1}.
 \end{equation}
Since the $(m-2)$-form on the right-hand side is closed and of rank at most two on $M$, from Lemma \ref{lemma:RhoRank2} it follows that $\psi_{i}$ is a Poisson structure of rank at most two, and hence is conformally invariant due to Lemma \ref{lemma:conformally}. This implies that we have Poisson structures on~$U$, of rank at most two, defined by
 \begin{equation}\label{eq:Pi-i}
 \pi_{i} := (-1)^{i-1}\,\alpha(X)\psi_{i}.
 \end{equation}
Here, $\alpha\in\Gamma(T^*U)$ is such that $\Omega = \alpha \wedge \beta$ on $U$, which exists by the independence of $h_1,\dots,h_{m-1}$. Now, since $\ii_X\dd h_j=0$ for all $j=1,\dots,m-1$, we have $\ii_{\pi_{i}^{\sharp}\dd h_{i}}\Omega = \dd h_{i}\wedge\ii_{\pi_{i}}\Omega = \alpha(X)\beta = \ii_{X}\Omega$. Hence, $\pi_{i}^{\sharp}\dd h_{i} = X$ and $\pi_{i}^{\sharp}\dd h_{j} = 0$ for $i \neq j$. Finally, taking into account that $\ii_{\pi_{j}}\ii_{\pi_{i}}\Omega = 0$ and $\ii_{\pi_{j}}\dd\ii_{\pi_{i}}\Omega = 0$, for all $i$ and $j$, the relations $[\pi_{i},\pi_{j}] = 0$ hold.
\end{proof}

Note that although the $\psi_{i}$s in \eqref{eq:Psi-i} are defined on the whole $M$, we only can assert by construction that the Poisson structures $\pi_i$ in \eqref{eq:Pi-i} exist on the same open set where the first integrals $h_1, \dots, h_{m-1}$ of $X$ are independent. However, in some cases, the Poisson structures $\pi_i$ are well defined on the whole manifold $M$. This occurs if, for instance, one can choose $\alpha$ in \eqref{eq:Pi-i} such that $\alpha(X)$ is constant, as we illustrate in the following examples.

\begin{Example}
On $\R{3}$, oriented with the Euclidean volume form $\Omega= \dd x \wedge \dd y \wedge \dd z$, consider the vector field governing the Euler's rigid body equations
 \begin{equation*}
 X = a yz \frac{\partial}{\partial x}+ b xz \frac{\partial}{\partial y} + c xy \frac{\partial}{\partial z},
 \end{equation*}
where $a$, $b$ and $c$ are nonzero constants related to the principal moments of inertia as follows: $a=\frac{I_2-I_3}{I_2I_3}$, $b=\frac{I_3-I_1}{I_1I_3}$ and $c=\frac{I_1-I_2}{I_1I_2}$. The vector field $X$ has the following first integrals, consisting of the energy and the square of angular moment:
 \begin{equation*}
 h_1 = \frac{1}{2}\bigg( \frac{x^2}{I_1} + \frac{y^2}{I_2} + \frac{z^2}{I_3} \bigg) \qquad \text{and} \qquad h_2 = \frac{1}{2}\big( x^2 + y^2 + z^2 \big).
 \end{equation*}
Since the first integrals $h_1$ and $h_2$ are independent on $U=\R{3} \setminus \{\text{coordinate axes}\}$, by Theo\-rem~\ref{teo:IntegrableHam}, the vector field $X$ is Hamiltonian on $U$ with respect to the commuting Poisson structures given by
 \begin{gather*}
 \pi_1 = x\frac{\partial}{\partial y}\!\wedge\! \frac{\partial}{\partial z} + y\frac{\partial}{\partial z}\!\wedge\! \frac{\partial}{\partial x}+z\frac{\partial}{\partial x}\!\wedge\! \frac{\partial}{\partial y},
 \qquad
 \pi_2 = \frac{x}{I_1}\frac{\partial}{\partial z}\!\wedge\! \frac{\partial}{\partial y} + \frac{y}{I_2}\frac{\partial}{\partial x}\!\wedge\! \frac{\partial}{\partial z}+\frac{z}{I_3}\frac{\partial}{\partial y}\!\wedge\! \frac{\partial}{\partial x}.
 \end{gather*}
However, note that these Poisson structures are well defined on the whole $\R{3}$, and hence the vector field $X$ is Hamiltonian on $\R{3}$.
\end{Example}

\begin{Example}\label{ex:TorusIntegrable}
Let $\mathbb{T}^2$ be the 2-torus with angular coordinates $(\varphi_1,\varphi_2)$ and volume form \mbox{$\Omega= \dd \varphi_1 \wedge \dd \varphi_2$}. For $m$ and $n$ coprime integers, consider the vector field
 \begin{equation*}
 X = \omega(n\varphi_1-m\varphi_2) \bigg( m\frac{\partial}{\partial \varphi_1} + n\frac{\partial}{\partial \varphi_2} \bigg),
 \end{equation*}
where $\omega \in C^\infty(\mathbb{R})$ is such that $\omega(t+2\pi) = \omega(t)$, for all $t\in \mathbb{R}$. For every $2\pi$-periodic $F \in C^\infty(\mathbb{R})$, we have that $h(\varphi_1,\varphi_2):= F(n\varphi_1-m\varphi_2)$ is a first integral of $X$. Furthermore, if $O \subset [0,2\pi]$ is the regular domain of $F$, then $U:=\{(\varphi_1,\varphi_2) \mid n\varphi_1-m\varphi_2\in O\}$ is the regular domain of $h$. By~Theorem \ref{teo:IntegrableHam}, the vector field $X$ is Hamiltonian with respect to the Poisson structure on $U$ given by
 \begin{equation*}
 \pi = \frac{\omega(n\varphi_1-m\varphi_2)}{F'(n\varphi_1-m\varphi_2)}\,\frac{\partial}{\partial \varphi_1}\wedge \frac{\partial}{\partial \varphi_2}.
 \end{equation*}
In principle, the Poisson structure $\pi$ is defined only on the proper open subset $U$ of $M$. But, if the function $\omega$ satisfies the condition $\int_0^{2\pi}\omega(t)\dd t=0$, then we can choose $F$ as a primitive of $\omega$, $F'=\omega$. In this case, the Poisson structure reduces to $\pi = \frac{\partial}{\partial \varphi_1}\wedge \frac{\partial}{\partial \varphi_2}$, which is well defined on the whole $\mathbb{T}^2$ even though the first integral $h(\varphi_1,\varphi_2)$ has singular points.
\end{Example}

We remark that relations \eqref{HamX} imply that each Poisson structure $\pi_i$ has $m-2$ Casimir functions given by the first integrals $h_1,\dots, h_{i-1},h_{i+1},\dots, h_{m-1}$ of $X$. Moreover, by \eqref{eq:Pi-i}, the rank of $\pi_i$ on $U$ is two except at the points where $X$ vanishes. Therefore, the \emph{set of singular points} of the Poisson structure $\pi_i$ on $U$ is the boundary of the open set in which $X$ is non-vanishing.

Now, we show that the bivector fields $\pi_1,\dots,\pi_{m-1}$ do not depend on the choice of the volume form. For any other volume form $\overline{\Omega}$ on $M$, there exists $\overline{\alpha} \in \Gamma(T^{\ast}U)$ such that $\overline{\Omega} = \overline{\alpha}\wedge\dd h_1\wedge\cdots\wedge\dd h_{m-1}$ on~$U$. Since $\alpha,\dd h_1,\dots,\dd h_{m-1}$ are independent on $U$, we have $\overline{\alpha} = f\alpha + g_1\dd h_1 + \cdots + g_{m-1}\dd h_{m-1}$, for some $f,g_1,\dots,g_{m-1} \in C^{\infty}(U)$. Moreover, it holds that $\overline{\Omega}=f\Omega$, and so $f$ is nowhere vanishing. By using $\overline{\alpha}$ and $\overline{\Omega}$, let us define $\overline{\psi}_{i}$ and $\overline{\pi}_{i}$ analogously as in~\eqref{eq:Psi-i} and \eqref{eq:Pi-i}, respectively. Consequently, $\overline{\psi}_{i}=\tfrac{1}{f}\psi_{i}$ and, taking into account that $\dd h_{i}(X)=0$, we~get
 \begin{equation*}
 \overline{\pi}_{i} = (-1)^{i-1}\,\overline{\alpha}(X)\overline{\psi}_{i} = (-1)^{i-1}f\,\alpha(X)\bigg( \frac{1}{f} \psi_{i} \bigg) = \pi_{i}.
 \end{equation*}

Finally, let us describe the dependence of $\pi_1,\dots,\pi_{m-1}$ on the choice of the first integrals $h_1,\dots,h_{m-1}$. If $\widetilde{h}_1,\dots,\widetilde{h}_{m-1}$ are (another) first integrals of $X$, independent on an open dense subset $V\subseteq M$, then $\dd \widetilde{h}_{i} = \sum_{j=1}^{m-1}a^{j}_{i}\dd h_{j}$, for some invertible matrix of functions $\big[a^{j}_{i}\big]$ on $U\cap V$. Hence, by straightforward computations, the Poisson structures $\widetilde{\pi}_i$ on $V$, defined from $\widetilde{h}_1, \dots, \allowbreak \widetilde{h}_{m-1}$ by \eqref{eq:Pi-i} and \eqref{eq:Psi-i}, satisfy
 \begin{equation}\label{eq:PiTilde}
 \widetilde{\pi}_i = \sum_{j=1}^{m-1}b^{i}_{j}\pi_j
 \end{equation}
on $U\cap V$, where $\big[b^{i}_{j}\big]$ is the inverse of the matrix $\big[a^{j}_{i}\big]$.

The above observations allows us to extend Theorem~\ref{teo:IntegrableHam} to not necessarily orientable manifolds.

\begin{Proposition}\label{prop:IntegrableHam}
Let $M$ be an $m$-dimensional manifold. Then, every vector field on $M$ admitting $m-1$ first integrals independent on a open dense subset $U\subseteq M$ is Hamiltonian with respect to pairwise commuting Poisson structures on $U$ of rank at most two.
\end{Proposition}
\begin{proof}
Let $X$ be a vector field on $M$ and $h_1, \dots, h_{m-1} \in C^{\infty}(M)$ first integrals independent on a open dense subset $U \subseteq M$. Let $\mathcal{U}=\{U_{\lambda}\}_{\lambda\in\Lambda}$ be a cover of $U$ by orientable open sets. For each $\lambda \in \Lambda$, fix a volume form $\Omega_{\lambda}$ on $U_{\lambda}$. By Theorem \ref{teo:IntegrableHam}, there exist commuting Poisson structures $\pi_{1}^{\lambda},\dots,\pi_{m-1}^{\lambda}$ on $U_{\lambda}$, of rank at most two, such that $\big(\pi_{i}^{\lambda}\big)^{\sharp}\dd h_{j} = \delta_{ij}X$, for $i,j=1,\dots,m-1$. Now, for fixed $\lambda,\lambda'\in\Lambda$, we have that $\pi_{i}^{\lambda}$ and $\pi_{i}^{\lambda'}$ agree on $U_{\lambda}\cap U_{\lambda'}$ since its constructions are independent of the choice of the volume forms $\Omega_{\lambda}$ and $\Omega_{\lambda'}$, respectively. Therefore, there exists $\pi_{i}\in\Gamma\big({\wedge}^2TM\big)$ such that $\pi_{i}|_{U_{\lambda}} = \pi_{i}^{\lambda}$, for all $\lambda \in \Lambda$. Finally, by construction, it follows that $\pi_1,\dots,\pi_{m-1}$ are commuting Poisson structures of rank at most two satisfying~\eqref{HamX}.
\end{proof}

We remark that if $M$ is non-orientable and connected and $X$ is non-trivial, then the Poisson structures obtained by the conclusion of Proposition \ref{prop:IntegrableHam} are singular. Indeed, their singular sets on $U$ are all the boundary of the open set where $X$ is non-vanishing, which is non-empty since~$M$ is non-orientable and $X$ admits the maximum number of first integrals.

{\bf Hamiltonization via invariant volume forms.} As an improvement of Theorem \ref{teo:IntegrableHam}, we proceed to show that we can relax the hypothesis on the number of first integrals: in the case when we have $m-2$ independent first integrals of a given vector field $X$, the existence of an $X$-invariant volume form is equivalent to the existence of a Poisson structure for which $X$ is an infinitesimal automorphism. Hence, the hamiltonization problem of $X$ turns into providing conditions for the triviality of its Poisson cohomology class.

\begin{Theorem}\label{teo:VolHam}
Let $X$ be a vector field on an orientable $m$-dimensional manifold $M$ admitting $m-2$ first integrals $c_1, \dots, c_{m-2} \in C^{\infty}(M)$ independent on an open dense subset $U \subseteq M$. Suppose that
 \begin{enumerate}\itemsep=0pt
 \item[$1)$] there exists an $X$-invariant volume form on $M$,
 \item[$2)$] the level sets of $\boldsymbol{c}:=(c_1,\dots,c_{m-2})\colon U \subseteq M \to \R{m-2}$ are connected and simply connected.
 \end{enumerate}
Then, the vector field $X$ is Hamiltonian on $U$ with respect to a unimodular Poisson structure on $M$ of rank two on $U$ and zero on $M \setminus U$, which admits $c_1,\dots,c_{m-2}$ as Casimir functions.
\end{Theorem}

In order to prove this theorem, we need the following:

\begin{Lemma}\label{lemma:1-1PoissVol}
Let $c_1,\dots,c_{m-2} \in C^{\infty}(M)$ be independent functions on an orientable manifold~$M$ of dimension $m$. Then, there exists a one-to-one correspondence between volume forms on~$M$ and unimodular Poisson structures of rank two on $M$ admitting $c_1,\dots,c_{m-2}$ as Casimir functions.
\end{Lemma}
\begin{proof}
Let us show that the correspondence between volume forms $\Omega$ on $M$ and Poisson structures $\pi$ of rank two such that each $c_i$ is a Casimir function is given by the relation
 \begin{equation}\label{eq:1-1PoissVol}
 \ii_{\pi}\Omega = \dd c_1 \wedge \cdots \wedge \dd c_{m-2}.
 \end{equation}
Indeed, the correspondence $\Omega \mapsto \pi$ follows from Lemmas \ref{lemma:RhoRank2} and \ref{lemma:HamRho}. Conversely, suppose we are given a Poisson structure $\pi$ of rank two with prescribed Casimir functions $c_1,\dots,c_{m-2}$ and fix a volume form $\Omega_0$. By Lemma \ref{lemma:HamRho}, we have $\ii_{\pi}\Omega_{0} \wedge \dd c_{i} = 0$, for all $i = 1,\dots,m-2$. Then, it follows from the independence of the $c_{i}$s that $\ii_{\pi}\Omega_{0} = f \dd c_1 \wedge \cdots \wedge \dd c_{m-2}$ for a certain nowhere vanishing $f \in C^{\infty}(M)$. By setting $\Omega := \tfrac{1}{f}\Omega_{0}$, the relation \eqref{eq:1-1PoissVol} follows.
\end{proof}

\begin{Remark}
The Poisson structures defined in \eqref{eq:1-1PoissVol} are called \emph{Flashcka--Ratiu Poisson structures}~\cite{Damianou-2012}.
\end{Remark}

\begin{proof}[Proof of Theorem \ref{teo:VolHam}]
Let $\Omega$ be a volume form such that ${\rm div}_{\Omega}(X)=0$. By Lemma \ref{lemma:1-1PoissVol}, the bivector field $\pi$ defined by \eqref{eq:1-1PoissVol} is a unimodular Poisson structure on~$M$ of rank two on~$U$ and zero on~$M \setminus U$, admitting $c_1,\dots,c_{m-2}$ as Casimir functions. So, the regular foliation of~$U$ consisting of the level sets of $\boldsymbol{c}\colon U \to \R{m-2}$ is the symplectic foliation of~$\pi$ on~$U$. By the $X$-invariance of $\Omega$, and since $c_1,\dots,c_{m-2}$ are first integrals of $X$, we get that $X$ is a tangent Poisson vector field of~$\pi$ on~$U$. Hence, the Poisson cohomology class of $X$ lies in the tangential Poisson cohomology in degree one, $[X] \in H^{1}_{{\rm tan}}(U,\pi) = {\mathfrak{poiss}_{{\rm tan}}(U,\pi)} \,\big/\, {\mathfrak{ham}(U,\pi)}$. Since the fibers of~$\boldsymbol{c}$ are connected and simply connected, we have $H^{1}_{{\rm tan}}(U,\pi)=0$, due to Proposition~\ref{prop:H1triv}. Therefore, the cohomology class of $X$ is trivial and $X\in\mathfrak{ham}(U,\pi)$.
\end{proof}

\begin{Remark}
The so-called \emph{Euler--Jacobi theorem} leads to another proof of Theorem \ref{teo:VolHam}. Indeed, that result implies that $X$ can be solved by quadratures along each fiber of $\boldsymbol{c}\colon U \subseteq M \to \R{m-2}$ (see, for example, \cite{BBM-11} and \cite[Theorem 1]{Kozlov-85}). Then, $X$ is locally Hamiltonian along each fiber of $\boldsymbol{c}|_{U}$, with respect to the symplectic structure induced by the restriction of $\pi$ in \eqref{eq:1-1PoissVol}. Since the \emph{fibers are simply connected}, the vector field $X$ is Hamiltonian along every whole fiber of $\boldsymbol{c}|_{U}$. Finally, since $\boldsymbol{c}\colon U\subseteq M \to \R{m-2}$ is a submersion with connected fibers, there exists a~Hamiltonian function $h$ for $X$ with respect to the Poisson structure $\pi$ \cite[Proposition~7.4]{DazHec-91}.
\end{Remark}

A remarkable difference between the Poisson structures obtained in Theorems \ref{teo:IntegrableHam} and \ref{teo:VolHam} is that the former encodes the zeroes of the vector field, while latter is always regular, so that the zeroes of the vector field are encoded in the Hamiltonian function. The regularity property in Theorem~\ref{teo:VolHam}, which is due to the independence of the first integrals, is in general a necessary condition, as the following example illustrates.

\begin{Example}\label{ex:IntSing}
Consider the functions $f,c\colon \R{3} \to \mathbb{R}$ given by $f= \tfrac{1}{2}\big(x^2+y^2+z^2\big)$ and $c = 2f^2$. The fibers of $c$ are connected and simply connected since they consist of 2-spheres and a single point. For the Euclidean volume form $\Omega$ on~$\R{3}$, the bivector field $\pi$ defined by $\ii_{\pi}\Omega = \dd c$ is a~Pois\-son structure with Casimir function $c$. Now, consider the vector field $X=-y\tfrac{\partial}{\partial x}+x\tfrac{\partial}{\partial y}$, for which $\Omega$ is invariant and $c$ is a first integral. Then, by Theorem \ref{teo:VolHam}, the vector field $X$ is Hamiltonian for $\pi$ on the open dense set $U:=\R{3} \setminus \{0\}$, where $\dd c$ is non-vanishing. However, \emph{the vector field $X$ is not Hamiltonian for $\pi$ on the whole $\R{3}$}.
In fact, no 1-form $\alpha$ satisfies $\pi^{\sharp}\alpha = X$ on~$\R{3}$. To see this, suppose such an $\alpha = \alpha_{1}\dd{x} + \alpha_{2}\dd{y} + \alpha_{3}\dd{z}$ exists. By restricting to $U$, we~get $\alpha_1|_{U} = (xz)/\big(2f^2\big) + \lambda xf$, $\alpha_2|_{U} = (yz)/\big(2f^2\big) + \lambda yf$ and $\alpha_3|_{U} = \big({-}y^2-x^2\big)/\big(2f^2\big) + \lambda zf$, for~some $\lambda \in C^{\infty}(U)$. The substitution $\lambda f = -z/\big(2f^2\big) + \mu$, for $\mu \in C^{\infty}(U)$, yields
 \begin{equation*}
 \alpha_1|_{U} = x\mu, \qquad \alpha_2|_{U} = y\mu, \qquad \alpha_3|_{U} = \frac{1}{f} + z\mu.
 \end{equation*}
From the first equation, it follows that $\alpha_{1}$ vanishes on the $yz$-plane, so there exists a global $\rho \in C^{\infty}\big(\R{3}\big)$ such that $\alpha_{1} = x\rho$. From here, we have $\rho|_{U} = \mu$. By the third equation, we get $\tfrac{1}{f} = (z\rho-\alpha_{3})|_{U}$, which is impossible, because ${1}/{f}$ cannot be smoothly extended to $\R{3}$.
\end{Example}

We remark that in the 3-dimensional case, Theorem \ref{teo:VolHam} recovers the hamiltonization criteria of \cite[first theorem]{Gao-2000}, where the topological hypotheses on the level sets defined by the first integrals were obviated.

Let us explain why the topological condition 2 for the fibers in Theorem \ref{teo:VolHam} is necessary in general. Suppose we are given a submersion $\boldsymbol{c}=(c_1, \dots, c_{m-2})\colon M \to \R{m-2}$ and a fiber-wise 1-form $\alpha$ along the fibers of $\boldsymbol{c}$ satisfying the following properties:
 \begin{enumerate}\itemsep=0pt
 \item[$(a)$] $\alpha$ is nowhere vanishing;
 \item[$(b)$] $\alpha$ is fiber-wise closed;
 \item[$(c)$] for every nowhere vanishing $f\in C^{\infty}(M)$, the 1-form $f\alpha$ is not fiber-wise exact.
 \end{enumerate}
In this setting, for any volume form $\Omega$, the vector field $X$ defined by $\ii_{X}\Omega = \alpha \wedge \dd c_{1}\wedge \dots \wedge \dd c_{m-2}$, has $\Omega$ as invariant volume form, but $X$ cannot be Hamiltonian for a Poisson structure $\pi$ with Casimir functions $c_1, \dots, c_{m-2}$. In fact, by Lemma~\ref{lemma:1-1PoissVol}, any two of such Poisson structures differ by multiplication of nowhere vanishing functions. Therefore, if $\beta$ is a fiber-wise 1-form such that $\pi^{\sharp}\beta = X$, then $\beta = f\alpha$ for some nowhere vanishing $f \in C^{\infty}(M)$. So, by property~($c$) on $\alpha$, the vector field $X$ is not Hamiltonian for $\pi$. We illustrate this situation in dimension three with a~couple of examples.

\begin{Example}\label{exm:37}
Consider the vector field $X = \tfrac{\partial}{\partial z}$ on $M = \R{3} \setminus\{z\text{-axis}\}$. Clearly, $c= \frac{1}{2}\big(x^2 + y^2\big)$ is a first integral for $X$ and the Euclidean volume form on $\R{3}$ is $X$-invariant. Observe that every nowhere vanishing Poisson structure $\pi$ with Casimir function $c$ is of the form
 \begin{equation*}
 \pi = \frac{1}{f}\bigg(x\frac{\partial}{\partial y}- y\frac{\partial}{\partial x}\bigg)\wedge\frac{\partial}{\partial z},
 \end{equation*}
for a nowhere vanishing function $f \in C^{\infty}(M)$. Moreover, the vector field $X$ is an infinitesimal Poisson automorphism for $\pi$ if and only if $\frac{\partial f}{\partial z} = 0$, that is, $f = f(x,y)$. Note that a 1-form $\alpha \in \Gamma(T^{*}M)$ satisfies $\pi^{\sharp}\alpha = X$ if and only if
 \begin{equation*}
 \alpha = f\,\frac{x\dd y - y\dd x}{x^2+y^2} \,+\, g\dd c,
 \end{equation*}
for some $g \in C^{\infty}(M)$. We claim that there is no $\alpha$ of the above form which is exact. Indeed, let $\gamma$ be the unitary circle on the $xy$-plane centered at the origin. Then, $\int_{\gamma}\alpha = \int_{\gamma}f \neq 0$ since~$f$ is nowhere vanishing and $\gamma$ is connected. In particular, $\alpha$ is not exact. Therefore, $X$ is not a~Hamiltonian vector field for Poisson structures with Casimir $c$.
\end{Example}

Although the previous example exhibits the necessity of condition 2 in Theorem \ref{teo:VolHam}, the vector field $\tfrac{\partial}{\partial z}$ obviously admits a maximal number of independent first integrals and hence is hamiltonizable. In the following example, we present a vector field that admits an invariant volume form, but only a first integral.

\begin{Example}\label{ex:Torus1}
Fix $\lambda \in \mathbb{R}$ and set $F:= x^2+y^2+z^2+1$. On $\R{3}$, the vector field
 \begin{equation}\label{EcXlambda}
 X_{\lambda} = (2xz+\lambda y)\frac{\partial}{\partial x} + (2yz-\lambda x)\frac{\partial}{\partial y} + \big(1-x^2-y^2+z^2\big)\frac{\partial}{\partial z}
 \end{equation}
has the first integral $c = \big(x^2\!+\!y^2\big)/F^2$ and admits the invariant volume form $\Omega = \frac{1}{F^2}\,\dd x \!\wedge \!\dd y \!\wedge\! \dd z$. Moreover, the level surface $T_{a} := c^{-1}(a)$ is a torus, for $a > 1/4$. Note that if $\lambda$ is rational, then the orbits of $X_{\lambda}$ along $T_{a}$ are closed paths and consequently the vector field $X_{\lambda}$ admits a~second first integral independent with $c$. By Theorem \ref{teo:IntegrableHam}, $X_\lambda$ is Hamiltonizable. If, instead, $\lambda$~is irrational, then the orbits of $X_{\lambda}$ are dense curves on $T_{a}$. This implies that any first integral of $X_{\lambda}$ is constant on $T_{a}$ and hence cannot be independent with $c$. In particular, it cannot be hamiltonized by using Theorem \ref{teo:VolHam}.
\end{Example}

\begin{Remark}
Even though the vector field $X_{\lambda}$ in \eqref{EcXlambda} does not admit a second first integral for $\lambda \in \mathbb{R} \setminus \mathbb{Q}$, there exists a Poisson structure with respect to which $X_{\lambda}$ is Hamiltonian (see Example \ref{ex:Torus2}).
\end{Remark}

Now, by using Theorem \ref{teo:VolHam}, we present families of linear vector fields admitting hamiltonization. Moreover, we give explicit formulas for the Hamiltonian function, the Poisson structure and the leaf-wise symplectic form based on the following fact: if $\Pi=\frac{1}{2}\Pi^{ij}\,\partial/\partial{x^{i}} \wedge \partial/\partial{x^{j}}$ is a~Poisson structure of \emph{rank at most two}, then the corresponding leaf-wise symplectic form is given by
 \begin{equation}\label{eq:FormaSimplect2}
 \omega_{S} =
 \begin{cases}
 -\dfrac{1\ }{|\Pi|^{2}} \Pi^{ij}\dd{x^{i}} \wedge \dd{x^{j}}\,\big|_{S} & \text{if} \  \dim{S} = 2, \\
 0 & \text{if} \ \dim{S} = 0,
 \end{cases}
 \end{equation}
where $|\Pi|^{2} := \sum_{1 \leq i < j \leq m}\big(\Pi^{ij}\big)^{2}$ and $S$ is a symplectic leaf of~$\Pi$.

\begin{Example}
Let $\R{n} = \{x = (x_{1}, \dots, x_{n})\}$ be the Euclidean vector space.
\begin{enumerate}\itemsep=0pt
 \item Consider the linear vector field $X(x) := Ax \cdot \frac{\partial}{\partial x}$ on $\R{n}$, associated with an $n \times n$ mat\-rix~$A$ such that $\operatorname{tr}A = 0$ and $\operatorname{rank}A \leq 2$. By the first condition, the canonical volume form on $\R{n}$ is $X$-invariant. The second one implies that there exist (linear) independent first integrals $c_1(x) := v_1 \cdot x, \dots, c_{n-2}(x) := v_{n-2} \cdot x$ of $X$, associated with some independent $v_1, \dots, v_{n-2} \in \ker{A^{\top}}$. Hence, by Theorem \ref{teo:VolHam}, the vector field $X$ is hamiltonizable on the whole $\R{n}$. Furthermore, by \eqref{eq:FormaSimplect2}, the Poisson structure in \eqref{eq:1-1PoissVol} is constant and given by
 \begin{equation}\label{EcPpoisson}
 \pi = \sum_{1 \leq i < j \leq n} (-1)^{i+j}\,\det{P_{[i,j]}} \frac{\partial}{\partial{x^{i}}} \wedge \frac{\partial}{\partial{x^{j}}},
 \end{equation}
 where $P_{[i,j]}$ denotes the $(n-2) \times (n-2)$ submatrix of $P=(v_{1} \cdots v_{n-2})$ without the rows~$i$ and~$j$. The symplectic foliation consists of the 2-dimensional planes given as the common level sets of the $c_i$s equipped with the constant symplectic structure induced by the restriction of $\omega = \omega_{ij}\dd{x^{i}} \wedge \dd{x^{j}}$, where
 \begin{equation}\label{EcPformaSym}
 \omega_{ij} = (-1)^{i+j+1} \frac{\det{P_{[i,j]}}}{|\pi|^{2}}, \qquad
 |\pi|^{2} := \!\!\sum_{1 \leq i < j \leq n}\!\!\big(\det{P_{[i,j]}}\big)^{2}, \qquad
 1 \leq i < j \leq n.
 \end{equation}
 In particular, a Hamiltonian function for $X$ is given by $h(x) = \frac{1}{2}x^{\top}(WA)x$, where $W=[\omega_{ij}]_{n \times n}$.

 \item Let $X (x) = Ax \cdot \frac{\partial}{\partial x}$ be a linear vector field on $\R{n}$, associated with a matrix $A$ such that $\operatorname{tr}A = 0$ and $\operatorname{rank}A = 3$. Suppose that the non-zero eigenvalues $\lambda_{1}$, $\lambda_{1}$, $\lambda_{3}$ of $A^{\top}$ are all real and distinct. Let $v_1$, $v_2$, $v_3$ be eigenvectors of $\lambda_1$, $\lambda_2$, $\lambda_3$, respectively; and $v_4,\dots,v_{n}$ a basis of $\ker{A}^{\top}$. As in above, the canonical volume form on $\R{n}$ is $X$-invariant and the linear functions $c_4(x) := v_4 \cdot x, \dots, c_{n}(x) := v_{n} \cdot x$ are independent first integrals of $X$. Moreover, the cubic function $c(x):=(v_1\cdot x)(v_2\cdot x)(v_3\cdot x)$ is also a first integral of $X$. This follows from ${\rm L}_{X}(v_i\cdot x) = \lambda_i v_i\cdot x$ and $\operatorname{tr}A = 0$. Hence, by Theorem \ref{teo:VolHam}, the vector field $X$ is hamiltonizable on the open set complementary to the union of the 2-planes generated by pairs of $v_1$, $v_2$, $v_3$. By \eqref{eq:FormaSimplect2}, the (quadratic) Poisson structure and the corresponding leaf-wise symplectic form are given by \eqref{EcPpoisson} and \eqref{EcPformaSym}, respectively, with $P = (u(x) \, v_{4} \,\cdots\, v_{n})$ and $u(x) = \sum_{{\rm cyclic}} (v_1 \cdot x)(v_2 \cdot x)\,v_3$. A Hamiltonian function for $X$ is given analogously as in above.

 \item Fix $a,v \in \R{3}$ such that $\|a\| > \|v\|$ and $a \cdot v=0$. Consider the linear vector field on $\R{4}=\big\{(x,y) \mid x \in \R{3}, y \in \mathbb{R}\big\}$ given by
 \begin{equation*}
 X(x,y) := (a \times x + yv)\cdot\frac{\partial}{\partial x}+(v\cdot x)\frac{\partial}{\partial y}.
 \end{equation*}
 We show that $X$ can be hamiltonized in two different ways by using Theorem \ref{teo:VolHam}. First, note that the quadratic function $q\colon \R{4} \setminus \{0\} \to \mathbb{R}_{s}$, $q(x,y) := \tfrac{1}{2}\big(x\cdot x - y^2\big)$, is a first integral of $X$ whose level sets $q^{-1}(s)$ are diffeomorphic to $\mathbb{S}^2 \times \mathbb{R}$ if $s \geq 0$ and to $\mathbb{S}^{0} \times \R{3}$ if $s < 0$. Now, the hypotheses on $a$ and $v$ imply that there exists $w \in \R{3}$, with $\|w\| < 1$, such that $w \times a = v$. Thus, the linear function $\ell(x,y) := w \cdot x - y$ is also a first integral of $X$. Finally, set $U:=\R{4} \setminus \{(x,y) \mid x - yw = 0 \}$. It follows from the property $\|w\| < 1$ that the level sets of $\boldsymbol{c}:=(q,\ell)\colon U \to \R{2}$ are diffeomorphic to $\mathbb{S}^2$, which is connected and simply connected. Since the Euclidean volume on $\R{4}$ is $X$-invariant, Theorem~\ref{teo:VolHam} implies that~$X$ is a~Hamiltonian vector field on $U$ with respect to a linear Poisson structure. On~the other hand, observe that the linear function $\varphi(x,y) := a\cdot x$ is also a first integral of~$X$. By restricting to the level sets of $(\ell,\varphi)\colon \R{4} \to \R{2}$, we get from the first incise that~$X$ is Hamiltonian with respect to a constant Poisson structure, with Hamiltonian function $h = \lambda\big(x\cdot x-y^2\big)$, for some constant factor $\lambda$ (see also \cite[Theorem 2]{Kozlov-92}).
 \end{enumerate}
\end{Example}

{\bf Hamiltonization via unimodularization.} A well-known fact about the hamiltonization problem, which can be seen as a trivial instance of Theorem~\ref{teo:VolHam}, is that a vector field~$X$ on a~simply connected 2-dimensional manifold is hamiltonizable if and only if it admits an invariant volume form $\Omega$ (we recall that every simply connected manifold is orientable). In this case, the Poisson structure~$\pi$ is defined by $\pi^{\sharp} = \big(\Omega^{\flat}\big)^{-1}$ and a Hamiltonian function is a primitive of $\ii_{X}\Omega$. Motivated by this situation, we formulate the following:

\begin{Theorem}\label{teo:UnimodHam}
Let $X$ be a vector field on an orientable $m$-dimensional manifold $M$. Suppose that there exists a volume form $\Omega$ such that
 \begin{equation}\label{EcXdRho}
 \ii_{X}\Omega = \dd\varrho,
 \end{equation}
for some $\varrho \in \Gamma\big({\wedge}^{m-2}T^{*}M\big)$. If $\varrho$ is of rank at most two and admits an integrating factor $a \in C^{\infty}(M)$, in the sense that $\dd(a\varrho) = 0$, then $X$ is Hamiltonian on the open set $U := \{a \neq 0\}$ with respect to a unimodular Poisson structure of rank at most two on $M$.
\end{Theorem}
\begin{proof}
First, by Lemma \ref{lemma:RhoRank2}, the bivector field $\pi$ defined by $\ii_{\pi}\Omega = a\varrho$ is a unimodular Poisson structure on $M$ of rank at most two. Now, from the integrating factor property, the function $h = {1}/{a}$ satisfies $\ii_{X}\Omega = \dd{h} \wedge (a\varrho)$. Hence, by Lemma \ref{lemma:HamRho}, the vector field $X$ is Hamiltonian for~$\pi$ on $U$, with Hamiltonian function $h$.
\end{proof}

\begin{Remark}
Observe that Theorem \ref{teo:UnimodHam} is a hamiltonization criteria via ``unimodularization'': the Poisson structure $\pi$ is unimodular and $X$ is its \emph{modular vector field} with respect to~$\Omega$~\cite{We-97}. Some relations between Hamiltonian vector fields with invariant volume forms and the unimodularity of the corresponding (almost-)Poisson structure in the non-holonomic case are given in \cite{FGM-15}.
\end{Remark}

Note that in the particular case when the integrating factor $a$ is nonwhere vanishing, then~$X$ is Hamiltonian on the whole $M$. On the other hand, a necessary condition for \eqref{EcXdRho} is that ${\rm div}_{\Omega}{X} = 0$ holds on $M$. This is sufficient if the de Rham cohomology of $M$ in degree $m-1$ is trivial, for example, if $M$ is an sphere or a real projective space.

\begin{Corollary}
Let $M$ be a orientable $m$-dimensional manifold with $H^{m-1}_{{\rm dR}}(M) = 0$ and $X$ a~vector field on $M$ admitting an invariant volume form $\Omega$. Then, condition \eqref{EcXdRho} holds for some~$\varrho$. Furthermore, if $\varrho$ admits a nowhere vanishing integrating factor, then $X$ is hamiltonizable on~$M$ in the following cases:
 \begin{enumerate}\itemsep=0pt
 \item[$(a)$] The manifold $M$ is $3$-dimensional.
 \item[$(b)$] The manifold $M$ is $4$-dimensional and $\varrho \wedge \varrho = 0$.
 \end{enumerate}
\end{Corollary}

This follows from Theorem \ref{teo:UnimodHam} since the items ($a$) and ($b$) ensure that $\varrho$ has rank at most two.

\begin{Example}
Let $M$ be the 3-sphere or the 3-dimensional projective space. Then, for each divergence free vector field $X$ on $M$, with respect to a volume form $\Omega$ on $M$, there exists a~differential 1-form $\varrho$ on $M$ such that \eqref{EcXdRho} holds. If $\varrho$ admits a nowhere vanishing integrating factor, then $X$ is a Hamiltonian vector field on $M$ with respect to a unimodular Poisson structure.
\end{Example}

{\bf Hamiltonization via orientable foliations.} Now, we present a generalization of Theorem~\ref{teo:VolHam} in the framework of orientable foliations. Recall that a regular foliation $\mathcal{F}$ is said to be orientable if there exists a nowhere vanishing element $\eta \in \Gamma\big({\wedge}^{{\rm top}}T^{\ast}\mathcal{F}\big)$, called a leaf-wise volume form of~$\mathcal{F}$.

\begin{Theorem}\label{teo:DoubleFol}
Let $\mathcal{F}$ be an oriented regular foliation on $M$ of dimension $r$ and $X \in \Gamma(TM)$ a vector field tangent to $\mathcal{F}$ preserving the leaf-wise volume. Under the following two conditions, the vector field~$X$ is hamiltonizable by a unimodular Poisson structure of rank two on $M$:
 \begin{enumerate}\itemsep=0pt
 \item[$1.$] There exists a nowhere vanishing, leaf-wise closed and locally decomposable leaf-wise $(r-2)$-form $\beta \in \Gamma\big({\wedge}^{r-2}T^{*}\mathcal{F}\big)$ such that $\ii_{X}\beta = 0$.

 \item[$2.$] The foliated de Rham cohomology of the foliation $\mathcal{F}_{\beta}$ integrating $\ker\beta$ is trivial in degree one, $H^{1}_{{\rm dR}}(M,\mathcal{F}_{\beta}) = 0$.
 \end{enumerate}
\end{Theorem}

\begin{proof}
Let $\eta \in \Gamma(\wedge^{r}T^{*}\mathcal{F})$ be an $X$-invariant leaf-wise volume form on $\mathcal{F}$. Define the bivector field $\pi \in \Gamma\big({\wedge}^{2}T\mathcal{F}\big)$ by $\ii_{\pi}\eta := \beta$. Note that $\beta$ is of rank at most two since it is locally decomposable. By the fact that $\beta$ is leaf-wise closed, we get from Lemma \ref{lemma:RhoRank2} that $\pi$ is Poisson. Moreover, the symplectic foliation of $\pi$ is precisely $\mathcal{F}_{\beta}$. On the other hand, since $\ii_{X}\beta = 0$, we get that $X$ is tangent to $\mathcal{F}_{\beta}$ and that ${\rm L}_{X}\beta = 0$. From here and the $X$-invariance of $\eta$, we get that $X$ is an infinitesimal Poisson automorphism for $\pi$. Finally, taking into account that $H^{1}_{{\rm dR}}(M,\mathcal{F}_{\beta}) = 0$, we conclude from Proposition \ref{prop:PoissTan} that $X$ is Hamiltonian.
\end{proof}

\begin{Corollary}\label{cor:DoubleFol}
Let $X$ be a volume-preserving vector field on an oriented $m$-dimensional mani\-fold $M$. Then, the vector field $X$ is hamiltonizable by a unimodular Poisson structure of rank two on $M$ if:
 \begin{enumerate}\itemsep=0pt
 \item[$1.$] There exist independent $1$-forms $\alpha_1,\dots,\alpha_k\in\Gamma(T^{*}M)$ satisfying the integrability condition
 \begin{equation*}
 \dd\alpha_{i} \wedge \alpha_{1} \wedge \cdots \wedge \widehat{\alpha}_{i} \wedge \cdots \wedge \alpha_{k} = 0, \qquad \text{for all}\quad i=1,\dots,k,
 \end{equation*}
 and such that $X$ is tangent to the foliation $\mathcal{F}_{\alpha}$ integrating $\bigcap_{i=1}^{k}\ker\alpha_{i}$.
{\sloppy \item [$2.$] There exists a nowhere vanishing, leaf-wise closed and locally decomposable leaf-wise $(m-k-2)$-form $\beta \in \Gamma\big({\wedge}^{m-k-2}T^{*}\mathcal{F}_{\alpha}\big)$ such that $\ii_{X}\beta = 0$.

}

 \item[$3.$] The foliated de Rham cohomology of the foliation $\mathcal{F}_{\beta}$ integrating $\ker\beta$ is trivial in degree one, $H^{1}_{{\rm dR}}(M,\mathcal{F}_{\beta})=0$.
 \end{enumerate}
\end{Corollary}
\begin{proof}
For $\alpha := \alpha_1 \wedge \cdots \wedge \alpha_k$, the bivector field $\pi$ defined by $\ii_{\pi}\Omega = \alpha\wedge\beta$
is Poisson and $X \in \mathfrak{ham}(M,\pi)$.
\end{proof}

As a generalization of Theorem \ref{teo:VolHam}, the results presented in Theorem \ref{teo:DoubleFol} and Corollary~\ref{cor:DoubleFol} allow us to hamiltonize vector fields without having explicit first integrals (see, for ins\-tance, Theorem \ref{teo:RegFol} and subsequent comments below). Indeed, we can recover Theorem \ref{teo:VolHam} by setting~$\mathcal{F}$ as the trivial foliation and $\beta = \dd c_1 \wedge \cdots \wedge \dd c_{m-2}$, where $c_1, \dots, c_{m-2} \in C^{\infty}(M)$ are independent first integrals of the given vector field $X$.

\begin{Example}
Consider the manifold $M = N \times \big(\R{3} \setminus \{y_{3}\text{-axis}\}\big)$ endowed with the oriented regular foliation $(\mathcal{F}, \eta)$ given by
\[
\mathcal{F} := \bigcup_{x \in N}\{x\} \times \big(\R{3} \setminus \{y_{3}\text{-axis}\}\big), \qquad
\eta := \frac{{\rm e}^{a(x)}}{y_{1}^{2} + y_{2}^{2}}\, \dd{y_{1}} \wedge \dd{y_{2}} \wedge \dd{y_{3}}\,|_{\R{3} \setminus \{y_{3}\text{-axis}\}},
\]
and the $\mathcal{F}$-tangent and $\eta$-preserving vector field on $M$
 \begin{equation*}
 X = f_{1}(x)\, y_{3} \bigg( y_{1} \frac{\partial}{\partial y_{1}} + y_{2} \frac{\partial}{\partial y_{2}} \bigg) + f_{2}(x, y) \frac{\partial}{\partial y_{3}}.
 \end{equation*}
Here, $a, f_{1} \in C^{\infty}(N)$ and $f_{2} \in C^{\infty}(M)$ with $\frac{\partial f_{2}}{\partial y_{3}} = 0$. Note that the nowhere vanishing and leaf-wise closed differential 1-form
 \begin{equation*}
 \beta = {\rm e}^{b(x)}\,\frac{y_{2} \dd y_{1} - y_{1} \dd y_{2}}{y_{1}^2 + y_{2}^2}, \qquad b \in C^{\infty}(N),
 \end{equation*}
is such that $\ii_{X}\beta = 0$. Moreover, the kernel of $\beta$ integrates to the subfoliation $\mathcal{F}_{\beta}$ of $\mathcal{F}$ given by the fibers of the submersion $p\colon M \to N \times \mathbb{S}^{1}$, $\big( x; {y_{1}}/{\sqrt{y_{1}^{2} + y_{2}^{2}}}, {y_{2}}/{\sqrt{y_{1}^{2} + y_{2}^{2}}} \,\big)$.
Since the leaves of~$\mathcal{F}_{\beta}$ are connected and simply connected, we have $H^{1}_{{\rm dR}}(M, \mathcal{F}_{\beta}) = 0$ (see \cite[Proposition~7.4]{DazHec-91}). Therefore, by the proof of Theorem \ref{teo:DoubleFol}, the vector field~$X$ is Hamiltonian on $M$ with respect to the Poisson structure
 \begin{equation*}
 \pi_{\eta,\beta} = -{\rm e}^{b(x)}\bigg( y_{1}\frac{\partial}{\partial y_{1}} + y_{2}\frac{\partial}{\partial y_{2}} \bigg) \wedge \frac{\partial}{\partial y_{3}}.
 \end{equation*}
Furthermore, in the case when $f_{2}$ depends radially on the variables $y_1$, $y_2$, that is, $f_{2}(x,y) = g\big(x; y_{1}^{2} + y_{2}^{2}\big)$ for some $g\in C^{\infty}(N\times\mathbb{R})$, a Hamiltonian function for $X$ is given by $h(x,y) = -\frac{1}{2}{\rm e}^{-b(x)}\big[G\big(x,\ln\big(y_1^2+y_2^2\big)\big) + y_{3}^{2}f_{1}(x)\big]$, where $G\in C^{\infty}(N\times\mathbb{R}_{t})$ is such that ${\partial G}/{\partial t} = g\big(x,{\rm e}^{t}\big)$.
\end{Example}

{\bf The case of integrable vector fields in the broad sense.} A vector field $X$ on an $m$-dimensional manifold $M$ is said to be \emph{integrable in the broad sense} \cite[Definition 1]{Bog98} if, for some $0\leq k\leq m-1$, it admits $k$ functionally independent first integrals and a $(m-k)$-dimensional abelian Lie algebra $\mathfrak{g} \subseteq \Gamma(TM)$ of symmetries of $X$ that also preserve the given first integrals.

For $k=m-1$, a vector field $X$ is integrable in the broad sense if and only if it admits the maximum number of independent first integrals. Indeed, a $1$-dimensional Lie algebra of symmetries of $X$ is the one generated by $X$ itself. Now, for $k=m-2$, the integrability in the broad sense of $X$ is equivalent to the existence of $m-2$ independent first integrals and two independent commutative vector fields preserving $X$ and tangent to the level set of the given first integrals.

Observe that the hamiltonization criteria of Proposition \ref{prop:IntegrableHam} states that every vector field integrable in the broad sense for $k=m-1$ is hamiltonizable. Similarly, the hamiltonization criteria of Theorem \ref{teo:VolHam} can be applied for the case $k=m-2$.

\begin{Proposition}
Let $X$ be a vector field on an orientable $m$-dimensional manifold $M$ that is integrable in the broad sense for $k=m-2$. Then:
 \begin{enumerate}\itemsep=0pt
 \item[$1.$] There exists an $X$-invariant volume form.
 \item[$2.$] If the common level sets of the $m-2$ independent first integrals of $X$ are connected and simply connected, then $X$ is hamiltonizable on $M$.
 \end{enumerate}
\end{Proposition}
{\sloppy\begin{proof}
By hypothesis, there exist independent first integrals $c_1,\dots,c_{m-2}$ of $X$ and commutative vector fields $Y_1$ and $Y_2$ that are symmetries of $X$ and tangent to the level sets of $c_1,\dots,c_{m-2}$. Then, there exists a unique volume form $\Omega$ on $M$ such that $\ii_{Y_1}\ii_{Y_2}\Omega = \dd c_1\wedge\cdots\wedge\dd c_{m-2}$. Since~$Y_1$,~$Y_2$ and $c_1,\dots,c_{m-2}$ are all $X$-invariant, so it is $\Omega$, which proves the item 1. From here and Theorem \ref{teo:VolHam}, the item 2 follows.
\end{proof}}

\section{Hamiltonization through decomposable Poisson structures}\label{sec:HamGood}

In this section, we give necessary and sufficient conditions under which an arbitrary vector field~$X$ on a smooth manifold $M$ is hamiltonizable via a decomposable Poisson structure $\pi$, which is non-regular in general. Specifically, we look for Poisson structures of the form $\pi = Y \wedge X$, for a vector field $Y \in \Gamma(TM)$. Then, we apply these conditions to give some hamiltonization criteria, involving transversally invariant Riemannian metrics and submersions, to the case of proper actions of 1-dimensional Lie groups and for infinitesimal generators of proper actions of~general Lie groups.

Consider the open set $U \subseteq M$ where the vector field~$X$ is non-vanishing. Denote by $\mathcal{F}_{X}$ the 1-dimensional foliation of $U$ given by the orbits of $X$ and let $\nu_{X} := TU/T\mathcal{F}_{X}$ be its normal bundle. If $X$ is Hamiltonian with respect to a Poisson structure $\pi$ of rank at most two, then we have $\pi \wedge X \in \Gamma\big({\wedge}^{3}\big(\pi^{\sharp}(T^{*}M)\big)\big) = \{0\}$. This fact, and the nowhere vanishing property of $X$ on~$U$, imply that there exists a vector field $Y \in \Gamma(TU)$ such that $\pi = Y \wedge X$ on $U$. Moreover, it~is clear that $\pi|_{U}$ only depends of the class $s := [Y] \in \Gamma(\nu_{X})$ rather than on $Y$ itself. Finally, if~$h \in C^{\infty}(M)$ is a Hamiltonian of $X$, then $\dd h(Y) = 1$ on $U$.

The properties discussed in the previous paragraph motivate the notion of normal class relative to a vector field $X$, which encodes the freedom of the choice of the vector field $Y$ such that $\pi = Y \wedge X$.

\begin{Definition}
Let $X$ be a vector field on $M$. A normal class relative to $X$, or simply a \emph{normal class}, is an equivalence class of vector fields on $M$ with respect to the relation
 \begin{equation*}
 Y_1 \sim_{X} Y_2 \qquad \text{if and only if} \qquad (Y_2 - Y_1) \wedge X = 0.
 \end{equation*}
Moreover, a normal class $[Y]$ is said to be
 \begin{itemize}\itemsep=0pt
 \item \emph{invariant} if $[X,Y] \wedge X \wedge Y=0$,
 \item \emph{normalized} with respect to a first integral $h$ of $X$ (or, \emph{$h$-normalized}) if $\dd h(Y)X=X$.
 \end{itemize}
\end{Definition}

It is straightforward to verify that, for any normal class, the properties of invariance and of $h$-normalization are independent of the choice of a representative. So, we have the following reformulation for the hamiltonization problem through decomposable Poisson structures.

\begin{Theorem}\label{teo:NormalSection}
A vector field $X$ is \emph{hamiltonizable} by a Poisson structure of the form $Y \wedge X$ on~$M$ if and only if $X$ admits an invariant normal class that is normalized with respect to a~first integral of $X$ that is regular on $\operatorname{supp}(X)$. In this case, given a first integral $h$ regular on $\operatorname{supp}(X)$, a Poisson structure $\pi$ on $M$ satisfying $\pi^{\sharp}\dd h=X$ is given by
 \begin{equation}\label{eq:NormalSection}
 \pi:=\begin{cases}
 Y_0\wedge X & \text{on} \  \operatorname{supp}(X),
 \\
 0 & \text{elsewhere},
 \end{cases}
 \end{equation}
where $Y_0$ is any vector field defined on a neighborhood of $\operatorname{supp}(X)$ satisfying $\dd h|_{\operatorname{supp}(X)}(Y_0)=1$ and $[X,Y_{0}] \wedge X = 0$.
\end{Theorem}

In other words, the hamiltonization problem for $X$ by decomposable Poisson structures $Y \wedge X$ translates into the construction of normalized invariant normal classes.

Note that for every Poisson structure of the form $\pi = Y \wedge X$, for which $X$ is Hamiltonian, it~follows from \eqref{eq:NormalSection} that the zeroes of $X$ and $\pi$ agree. Consequently,
 \begin{equation*}
 \operatorname{rank}{\pi_{p}} =
\begin{cases}
 2 & \text{if} \  X_{p} \neq 0,
 \\
 0 & \text{if} \  X_{p} = 0,
 \end{cases}
 \end{equation*}
which implies that the set of singular points of $\pi=Y\wedge X$ is the boundary of the open set in~which $X$ is non-vanishing. In particular, the Poisson structure $\pi$ is regular if and only if, on~each connected component of $M$, the vector field $X$ is nowhere-vanishing or zero.

We have divided the proof of Theorem \ref{teo:NormalSection} into a sequence of lemmas. As a first step, observe that one can associate to each normal class $s = [Y]$ relative to a vector field $X$ on $M$ a well-defined bivector field $\pi_{s}$ given by $\pi_{s} := Y \wedge X$.

\begin{Lemma}\label{lemma:OneToOne}
The assignment $s \mapsto \pi_{s}$ is a one-to-one correspondence between normal classes relative to $X$ and bivector fields decomposable by $X$. Moreover,
 \begin{enumerate}\itemsep=0pt
 \item[$(i)$] the class $s$ is invariant if and only if $\pi_{s}$ is Poisson,
 \item[$(ii)$] the class $s$ is $h$-normalized if and only if $\pi_{s}^{\sharp}\dd h = X$.
 \end{enumerate}
\end{Lemma}
\begin{proof}
By definition of normal class, the assignment $s \mapsto \pi_{s}$ is injective. For the surjectivity, simply note that the decomposable bivector $\pi = Y \wedge X$ is the image of the normal class defined by $Y$. On the other hand, if $Y$ is a representative of the normal class $s$, then $[\pi_{s},\pi_{s}] = 2 [X,Y] \wedge X \wedge Y$, so the Poisson property for $\pi_{s}$ is equivalent to the invariance of $s$. Finally, if $h$ is a first integral of $X$, then we have $\pi_{s}^{\sharp}\dd h = \dd h(Y) X$, which implies that $s$ is $h$-normalized if and only if $\pi_{s}^{\sharp}\dd h = X$.
\end{proof}

Now, we show that the first integrals normalizing some normal class of a vector field $X$ are those that are regular on the support of $X$.

\begin{Lemma}\label{lemma:NormalGood}
Let $X$ be a vector field with first integral $h$. Then, an $h$-normalized class relative to $X$ exists if and only if $h$ is regular on $\operatorname{supp}(X)$. In this case, if $U$ is the open set where $X$ is non-vanishing, then:
 \begin{enumerate}\itemsep=0pt
 \item[$1.$] Every vector field $Y_0$ defined on a neighborhood of $\operatorname{supp}(X)$ satisfying $\dd h|_{U}(Y_0) = 1$ determines a unique $h$-normalized normal class $s_{0}$ on $M$.
 \item[$2.$] If, in addition, the identity $[X,Y_0] \wedge X \wedge Y_0 = 0$ holds on $U$, then $s_{0}$ is invariant on $M$.
 \end{enumerate}
\end{Lemma}
\begin{proof}
If $Y$ represents an $h$-normalized normal class, then $\dd h(Y) = 1$ on $U$, which implies that $\dd h(Y) = 1$ holds at $\overline{U} = \operatorname{supp}(X)$. Therefore, $\dd h$ does not vanish on $\overline{U}$. Now, suppose that $h$ is regular on $\overline{U}$ and let $Y_0 \in \Gamma(TV)$ be a vector field defined on a neighborhood $V$ of $\overline{U}$ such that $\dd h|_{U}(Y_0) = 1$. Since $\overline{U} \subseteq V$, we can find disjoint open sets $W$ and $W'$ satisfying $\overline{U} \subseteq W$ and $M \setminus V \subseteq W'$, together with a smooth function $\mu \in C^{\infty}(M)$ such that $\mu|_{\overline{U}} = 1$ and $\mu|_{M \setminus W} = 0$. Then, a smooth vector field $Y$ on $M$ is well defined by $Y := \mu Y_0$ on~$V$ and $Y := 0$ on~$W'$. Since~$Y$ agrees with $Y_0$ on $U$, we have $\dd h(Y)X = X$. In other words, the normal class $s_0:=[Y]$ is $h$-normalized.
Finally, since $X(\mu) = 0$ and $X$ vanishes outside $U$, the invariance of $s_0$ follows from $[X,Y_0] \wedge X \wedge Y_0 = 0$ on $U$.
\end{proof}

Finally, we describe a reformulation of the invariance property for normalized classes.

\begin{Lemma}\label{lemma:NormalInv}
Let $X$ be a vector field on $M$ with first integral $h$. Let also $s$ be an $h$-normalized normal class and $Y$ a representative of $s$. Then, the following assertions are equivalent:
 \begin{enumerate}\itemsep=0pt
 \item[$1.$] The class $s$ is invariant.
 \item[$2.$] The identity $[X,Y]\wedge X = 0$ holds on $M$.
 \item[$3.$] There exists $a \in C^{\infty}(U)$ such that $[X,Y] = aX$ holds on $U$.
 \end{enumerate}
Here, we denote by $U$ the open subset of $M$ where $X$ is non-vanishing.
\end{Lemma}
\begin{proof}
Since, by definition, $X$ is nowhere vanishing on $U$, it readily follows that 2 and 3 are equivalent. Moreover, it is clear that 2 implies 1. So, it is left to show that 1 implies 2. Observe that $0 = \ii_{\dd h}([X,Y]\wedge X\wedge Y) = \dd h([X,Y]) X\wedge Y - \dd h(X)[X,Y]\wedge Y + [X,Y]\wedge(\dd h(Y)X) = \dd h([X,Y]) X\wedge Y + [X,Y]\wedge X$, where in the last step, we have applied that $\dd h(X) = 0$ and $\dd h(Y)X = X$. Since $\dd h([X,Y]) X \wedge Y = 0$ automatically holds on $M \setminus U$, it suffices to show that $\dd h([X,Y]) = 0$ on $U$. From $\dd h(Y)X = X$, we get that $\dd h(Y) = 1$ on $U$. Moreover, by the Koszul's formula, $\dd h([X,Y]) = {\rm L}_{X}(\dd h(Y)) - {\rm L}_{Y}(\dd h(X)) - \dd^2h(X,Y)$. Since $\dd h(X)$ and $\dd h(Y)$ are constant on $U$, we conclude that $\dd h([X,Y]) = 0$ on $U$, as desired.
\end{proof}

\begin{proof}[Proof of Theorem \ref{teo:NormalSection}]
The first assertion follows from Lemmas \ref{lemma:OneToOne} and \ref{lemma:NormalGood}. So, it remains to show that $\pi$ in \eqref{eq:NormalSection} defines a Poisson structure on $M$ satisfying $\pi^{\sharp}\dd h=X$. To see this, note that $\pi$ agrees on $M$ with the \emph{smooth} bivector field $\pi_{s_0}:=Y\wedge X$, where $Y\in\Gamma(TM)$ and $s_0$ are given as in the proof of Lemma \ref{lemma:NormalGood}. The fact that $\pi_{s_0}$ is Poisson and satisfies $\pi^{\sharp}_{s_0}\dd h=X$ follows from Lemmas \ref{lemma:OneToOne} and \ref{lemma:NormalInv}.
\end{proof}

\begin{Remark}
There is a cohomological approach to the proof of Lemma \ref{lemma:NormalInv}, where the non-trivial part is to show that $[X,Y] \wedge X \wedge Y = 0$ implies $[X,Y] \wedge X = 0$. Indeed, set $\pi := Y \wedge X$. Since $\dd h(Y)X = X$ and $[X,Y] \wedge X \wedge Y = 0$, we get from Lemma~\ref{lemma:OneToOne} that $\pi$ is Poisson and~$X$ is Hamiltonian. In particular, $X$ is an infinitesimal Poisson automorphism for $\pi$, which leads to $[X,Y] \wedge X = 0$.
\end{Remark}

Now, following \cite{LPS-12}, we use Theorem \ref{teo:NormalSection} to formulate a hamiltonization criteria of vector fields on orientable manifolds, that are \emph{a priori} not endowed with a first integral.

\begin{Theorem}\label{teo:DiverInt}
Let $X$ be a vector field on $M$ admitting an invariant volume form $\Omega$ and a vector field $Z$ satisfying $[X,Z]=\lambda X$, for some $\lambda \in C^{\infty}(M)$. Then, the function $h:=\operatorname{div}_{\Omega}(Z)-\lambda$ is a first integral of $X$. Moreover, if ${\rm L}_{Z}h$ is non-zero on $\operatorname{supp}(X)$, then $X$ is Hamiltonian on $M$ with respect to the Poisson structure of rank at most two defined by
 \begin{equation}\label{eq:PiDiverInt}
 \pi :=
 \begin{cases}
 \dfrac{1}{{\rm L}_{Z}h} Z \wedge X & \text{on} \  \operatorname{supp}(X),
 \\
 0 & elsewhere.
 \end{cases}
 \end{equation}
\end{Theorem}
\begin{proof}
First, we show that $h$ is a first integral of $X$. Indeed, since $\Omega$ is $X$-invariant and~$Z$ is a symmetry of $X$, we have ${\rm L}_{X}({\rm div}_{\Omega}(Z)\Omega) = {\rm L}_{[X,Z]}\Omega + {\rm L}_{Z}({\rm L}_{X}\Omega) = \dd \lambda\wedge\ii_{X}\Omega + \lambda\dd \ii_{X}\Omega = ({\rm L}_{X}\lambda)\Omega + \lambda({\rm L}_{X}\Omega) = {\rm L}_{X}(\lambda\Omega)$, which implies that ${\rm L}_{X}(h\Omega)=0$. Moreover, we also have $({\rm L}_{X}h)\Omega = {\rm L}_{X}(h\Omega) - h\,{\rm L}_{X}\Omega = 0$, showing that $h$ is a first integral of $X$.
Now, consider the open set $V := \{{\rm L}_{Z}h\neq 0\}$ and the vector field $Y_0:=\frac{1}{{\rm L}_{Z}h}Z$ on $V$. Then, $\dd h|_{V}(Y_0)=1$ and $X \wedge [X,Y_0] \wedge Y_0 = -\tfrac{\dd h([X,Z])}{{\rm L}_{Z}h^3}X \wedge(\lambda X)\wedge Z = 0$. So, by Lemma \ref{lemma:NormalGood}, the vector field $Y_0$ defines a unique normal class for $X$, which is $h$-normalized and invariant. Hence, by Theorem \ref{teo:NormalSection}, the vector field $X$ is Hamiltonian on $M$ with respect to the Poisson structure in \eqref{eq:PiDiverInt}.
\end{proof}

We left for Section \ref{subsec:Metric} the discussion of further applications of Theorem \ref{teo:NormalSection} to the problem of hamiltonization by means of possibly non-regular first integrals.

{\bf The case of regular first integrals.} In this part we review under the light of Theo\-rem~\ref{teo:NormalSection} some known hamiltonization criteria that require the existence of \emph{regular} first integrals, that is, without critical points. This requirement, however, imposes strong topological conditions on the manifold: it excludes, for instance, a global hamiltonization of vector fields on compact manifolds.

As a consequence of Theorem \ref{teo:NormalSection}, we recover the following hamiltonization criteria \cite{AHA-06,Hojman-96}, which implicitly require the regularity of the first integral.

\begin{Corollary}\label{cor:Hojman}
Let $X$ be a vector field admitting a regular first integral $h$.
 \begin{enumerate}\itemsep=0pt
 \item[$(a)$] If a vector field $Y$ is such that $\dd h(Y) = 1$ and $[X,Y] = aX$, for some $a \in C^{\infty}(M)$, then $\pi := Y \wedge X$ is Poisson and satisfies $\pi^{\sharp}\dd h = X$.
 \item[$(b)$] If a vector field $Z$ is such that $\dd h(Z)$ is nowhere vanishing and $[X,Z] = pX + qZ$, for some $p,q \in C^{\infty}(M)$, then $\pi := \tfrac{1}{\dd h(Z)}Z \wedge X$ is Poisson and satisfies $\pi^{\sharp}\dd h = X$.
 \end{enumerate}
\end{Corollary}
\begin{proof}
For item ($b$), it is straightforward to verify that the normal class of $Y := \tfrac{1}{\dd h(Z)}Z$ is invariant and $h$-normalized. Hence, by Theorem \ref{teo:NormalSection}, the bivector field $\pi := Y \wedge X = \tfrac{1}{\dd h(Z)}Z\wedge X$ is Poisson and satisfies $\pi^{\sharp}\dd h = X$. Moreover, item ($a$) readily follows from ($b$).
\end{proof}

The hamiltonization criteria in the item ($a$) of Corollary~\ref{cor:Hojman}, which is known as the \emph{Hojman construction}, appears for the first time in~\cite{Hojman-96}. An intrinsic formulation of ($a$) is found in~\cite{AHA-06}, where statement ($b$) is presented as a generalization of~($a$). However, we remark that ($b$) \emph{is not} a generalization but a reformulation of ($a$), since it does not provide new solutions to the hamiltonization problem at all: by Lemma~\ref{lemma:NormalInv}, every vector field $X$ satisfying ($b$) also satisfies~($a$). Indeed, by taking $Y := \tfrac{1}{\dd h(Z)}Z$, we get $\dd h(Y) = 1$ and $[X,Y] = aX$, for $a := {p}/{\dd h(Z)}$.

\begin{Example}
On $M=\{ x \in \R{n} \mid x_{i} \neq 0,\ i= 1,\dots,n\}$, with $n\geq2$, let us consider the vector field $X = \sum_{i=1}^{n} x_{i}F_{i}\frac{\partial}{\partial x_{i}}$, where $F_{1},\dots,F_{n}$ are homogeneous functions of degree $r$ on $M$ such that $F_{1} + \cdots + F_{n} = 0$. Note that the function $h(x)=x_{1} \cdots x_{n}$ is a regular first integral of~$X$ on~$M$. On the other hand, since $X$ is an homogeneous vector field of degree $r+1$, we have $[X,E] = -rX$, where $E = x_{1}\frac{\partial}{\partial x_{1}} + \cdots + x_{n}\frac{\partial}{\partial x_{n}}$ is the Euler vector field. Furthermore, the function $\dd h(E) = nh$ is nowhere vanishing on $M$. By Corollary~\ref{cor:Hojman}, the bivector field
 \begin{equation*}
 \pi = \frac{1}{nh} \sum_{1 \leq i < j \leq n} x_i x_j (F_j - F_i) \frac{\partial}{\partial x_i} \wedge \frac{\partial}{\partial x_j},
 \end{equation*}
is Poisson and satisfies $\pi^{\sharp}\dd h = X$. Therefore, $X$ is a Hamiltonian vector field on~$M$.
\end{Example}

\begin{Example}
Let $X = Ax\cdot\frac{\partial}{\partial x}$ be a linear vector field on $\R{n}_{x}$, associated with a real mat\-rix~$A$ and $E$ the Euler vector field on $\R{n}_{x}$. Observe that $[X,E] = 0$.
\begin{itemize}\itemsep=0pt
\item First, suppose that $A^{\top}$ admits two distinct real eigenvalues $\lambda_{1}$, $\lambda_{2}$, with eigenvectors~$v_{1}$,~$v_{2}$, respectively. Fix real numbers $r_{1}$, $r_{2}$ satisfying $\lambda_{1}r_{1} + \lambda_{2}r_{2} = 0$ and $r_{1} + r_{2} \neq 0$. Then, on the open set $U = \{x\in\R{n} \mid x\cdot v_{1}\neq0, \, x\cdot v_{2}\neq 0\}$, the function $h(x) = |v_{1} \cdot x|^{r_{1}} |v_{2}\cdot x|^{r_{2}}$ is a~regular first integral of $X$. Moreover, $\dd h(E)= (r_{1} + r_{2}) h$ is nowhere vanishing on~$U$. Then, by Corollary \ref{cor:Hojman}, the vector field $X$ is Hamiltonian on $U$ with respect to $\pi_{h} = \frac{1}{h}E \wedge X$, with Hamiltonian function $h$.
\item On the other hand, if the kernel of $A^{\top}$ is non-trivial, then for every non-zero $v \in \ker{A^{\top}}$ the linear function $f(x) = v \cdot x$ is a regular first integral of $X$. Moreover, $\dd f(E)= f$ is nowhere vanishing on $W = \{x\in\R{n} \mid x\cdot v\neq0\}$. By Corollary \ref{cor:Hojman}, the vector field $X$ is Hamiltonian on $W$ with respect to $\pi_{f} = \frac{1}{f}E \wedge X$, with Hamiltonian function $f$.
\end{itemize}
\end{Example}

Now, we present the following hamiltonization criteria for vector fields admitting a regular first integral and a suitable foliation.
This result is a generalization of \cite[Theorem 2]{FGS-05}, where it is formulated for nowhere vanishing vector fields with periodic flow.

\begin{Theorem}\label{teo:RegFol}
Let $X$ be a vector field on an $m$-dimensional manifold $M$ admitting a first integral $h\in C^{\infty}(M)$ with no critical points. Suppose that $X$ is tangent to a 2-dimensional foliation~$\mathcal{S}$ transversal to the level sets of $h$. Then, there exists a unique Poisson structure $\pi=Y \wedge X$, for some $Y \in \Gamma(TM)$, such that $\pi^{\sharp}\dd h=X$ on $M$, its symplectic foliation coincides with $\mathcal{S}$ on the open set $\{\operatorname{rank}{\pi} = 2\}$ and vanishing at the set $\{X=0\}$.
Moreover, on a~foliated coordinate chart $(U;\boldsymbol{x} = (x_1,\dots,x_{m-2}),\boldsymbol{y} = (y_1,y_2))$ for $\mathcal{S}$, such that $y_1=h$, we have $\pi|_{U}=\frac{\partial}{\partial y_1}\wedge X$.
\end{Theorem}
\begin{proof}
Since $h$ is transversal to $\mathcal{S}$ and has no critical points, around each point in $M$ there exist a foliated coordinate chart $(U,\boldsymbol{x},\boldsymbol{y})$ of $\mathcal{S}$ such that $y_1=h$ and the leaves of $\mathcal{S}|_{U}$ are the level sets of $\boldsymbol{x}$. By the fact that $X$ is tangent to $\mathcal{S}$, the functions $x_1,\dots,x_{m-2},h$ are independent first integrals of $X$. By Theorem \ref{teo:IntegrableHam}, the bivector field $\pi=\frac{\partial}{\partial y_1}\wedge X$ is the unique Poisson structure on $U$, of rank at most two, such that $\pi^{\sharp}\dd h=X$, with Casimir functions $x_{1}, \dots, x_{m-2}$. Now, for another foliated chart $(V,\widetilde{\boldsymbol{x}},\widetilde{\boldsymbol{y}})$ with $\widetilde{y}_{1}=h$,
we have that $\widetilde{x}_1,\dots, \widetilde{x}_{m-2},h$ are independent first integrals of $X$ on $V$. Then, there exists an invertible matrix of functions $A=[a^{j}_{i}]_{i,j=1,\dots,m-1}$ on $U\cap V$ satisfying $a^{j}_{m-1}=\delta^{j}_{m-1}$ and $\dd \widetilde{x}_{i} = \sum_{j=1}^{m-2}a^{j}_{i}\dd x_{j} + a^{m-1}_{i}\dd h$. Since the common level sets of $\widetilde{x}_{i},x_1,\dots,x_{m-2}$ are the (2-dimensional) leaves of $\mathcal{S}|_{U\cap V}$, it follows that $a^{m-1}_{i}\dd h\wedge\dd x_{1}\wedge\cdots\wedge\dd x_{m-2} = \dd \widetilde{x}_{i}\wedge\dd x_{1}\wedge\cdots\wedge\dd x_{m-2}=0$. This means that the matrix $A$ is of the form $\left(\begin{smallmatrix}
 \ast & 0 \\
 0 & 1
 \end{smallmatrix}\right)$,
with $\ast=\big[a^{j}_{i}\big]_{i,j=1,\dots,m-2}$. Consequently, its inverse matrix is of the same form. So, by \eqref{eq:PiTilde}, the Poisson structure $\widetilde{\pi}$ on $V$, obtained by Theorem \ref{teo:IntegrableHam} using $\widetilde{x}_1,\dots, \widetilde{x}_{m-2},h$, satisfies $\widetilde{\pi}|_{U\cap V}=\pi|_{U\cap V}$. This shows that $\pi$ is independent of the choice of the coordinate chart. Hence, the Poisson structure $\pi$ is global and satisfies $\pi^{\sharp}\dd h=X$. Furthermore, since $x_{1},\dots,x_{m-2}$ are Casimir functions of $\pi$, the symplectic foliation of $\pi$ coincides with $\mathcal{S}$ on the open set $\{\operatorname{rank}\pi = 2\}$. Also, by \eqref{eq:Pi-i}, the Poisson structure $\pi$ vanishes at $\{X=0\}$. Finally, by the partition of unity argument, there exists $Y \in \Gamma(TM)$ such that $\pi=Y \wedge X$.
\end{proof}

In the case when $X$ has no critical points, one can prove Theorem \ref{teo:RegFol} by using the Hojman construction: since the leaves of $\mathcal{S}$ are oriented by $X$ and ${\rm d}h|_{\mathcal{S}}$, there exists a nowhere vanishing vector field $Z\in\Gamma(T\mathcal{S})$ such that $X$ and $Z$ span $T\mathcal{S}$. Then, the involutivity of $\mathcal{S}$ implies that each condition of Corollary~\ref{cor:Hojman}(b) is satisfied. Therefore, $X$ is a Hamiltonian vector field with respect to $\pi = {1}/{\dd{h}(Z)}Z\ \wedge X$, that clearly satisfies the conclusion of Theorem~\ref{teo:RegFol}.

On the other hand, if $M$ is orientable, then the foliation $\mathcal{S}$ is given by the kernel of a differential $(m-2)$-form.
In this particular case, Theorem \ref{teo:RegFol} follows from Theorem \ref{teo:DoubleFol}.

 \subsection[Transversally invariant metrics, submersions and proper actions of 1-dimensional Lie groups]
 {Transversally invariant metrics, submersions and proper actions \\of 1-dimensional Lie groups}\label{subsec:Metric}

In this part, we develop some hamiltonization criteria in terms of transversally invariant Riemannian metrics (for a more general view of this notion and its geometric applications, see \cite{del_Hoyo_2018, PPT-2014}). Also, we consider the hamiltonization problem for vector fields that are vertical on the total space of a submersion, and we derive some results on hamiltonization of infinitesimal generators of proper actions that extend some of the results in \cite{FGS-05}, where similar approaches are applied in the case of periodic flow vector fields.

{\bf Hamiltonization via transversally invariant metrics.} Let $(M,g)$ be a Riemannian manifold and $X$ a vector field on $M$. First, we recall that the metric $g$ is said to be \emph{$X$-invariant} if ${\rm L}_{X}g = 0$. A weaker notion is that of transversal invariance: let $O$ be the orbit of $X$ through a point $m \in M$. Then, the metric $\eta \in \Gamma(T^{*}M \otimes T^{*}M)$ dual to $g$ induces a metric $\eta_{O}$ on the annihilator subbundle $TO^{\circ}$. We say that $g$ is \emph{transversally invariant} if, for each $m \in M$, the map $\big(\dd_{m}{\rm Fl}^{t}_{X}\big)^{*}\colon T_{{\rm Fl}^{t}_{X}(m)}O^{\circ}\to T_{m}O^{\circ}$ induced by the flow ${\rm Fl}^{t}_{X}$ of $X$ is an isometry, wherever is well defined. This is equivalent to require that the map $\overline{\dd_{m}{\rm Fl}^{t}_{X}}\colon T_{m}M / T_{m}O \to T_{{\rm Fl}^{t}_{X}(m)}M / T_{{\rm Fl}^{t}_{X}(m)}O$ is an isometry. Infinitesimally, this is just ${\rm L}_{X|_{O}}\eta_{O} = 0$ for every orbit $O$ of $X$. It is clear that every $X$-invariant metric is transversally invariant, but the converse is not true.

\begin{Lemma}\label{lemma:TransInv}
Let $X$ be a vector field and $g$ a transversally invariant Riemannian metric for~$X$. Then, for every first integral $h$ of $X$, its gradient $\nabla h:=\eta^{\sharp}(\dd h)$ commutes with $X$, and $g(\nabla h,\nabla h)$ is a first integral of $X$.
\end{Lemma}
\begin{proof}
Since $\dd h$ is $X$-invariant, we have $[X,\nabla h] = ({\rm L}_{X}\eta)^{\sharp}\dd h$. By the transversal invariance, we have $\beta([X,\nabla h])=0$ for each 1-form $\beta$ vanishing on $X$. Now, let $\theta$ be a 1-form orthogonal to those on the annihilator of $X$. Then, we have that ${\rm L}_{X}\theta$ is orthogonal to the annihilator of~$X$. In particular, $\theta$ and ${\rm L}_{X}\theta$ are orthogonal to $\dd h$. So, we have $\theta([X,\nabla h]) = {\rm L}_{X}(\eta(\dd h,\theta)) - \eta(\dd h,{\rm L}_{X}\theta) = 0$. This implies that $[X,\nabla h]=0$. Furthermore, we have ${\rm L}_{X}(g(\nabla h,\nabla h))={\rm L}_{X}\eta(\dd h,\dd h) + 2g({\rm L}_{X}\dd h,\dd h)=0$, where we have used the transversal invariance of $g$ and the invariance of $\dd h$.
\end{proof}

\begin{Theorem}\label{teo:TrInvM}
Let $X$ be a vector field on the manifold $M$ that admits
 \begin{enumerate}\itemsep=0pt
 \item[$1)$] a first integral $h \in C^{\infty}(M)$ that is regular on $\operatorname{supp}(X)$,
 \item[$2)$] a transversally invariant metric $g \in \Gamma(TM \otimes TM)$.
 \end{enumerate}
Then, the vector field $X$ is Hamiltonian on $M$ with respect to the Poisson structure of rank at most two defined by
 \begin{equation*}
 \pi :=
 \begin{cases}
 \dfrac{\nabla h}{g(\nabla h,\nabla h)} \wedge X & \text{on} \ \operatorname{supp}(X), \\
 0 & elsewhere.
 \end{cases}
 \end{equation*}
\end{Theorem}
\begin{proof}
Let $U$ and $V$ be the respective open sets where $X$ and $\dd h$ are nowhere vanishing. Define $Y_0 \in \Gamma(TV)$ by $Y_0 := \frac{\nabla h}{g(\nabla h,\nabla h)}$. Since $h$ is a first integral of $X$, we have $\dd h \in \Gamma(TO^{\circ})$ for each orbit $O$ of $X$. From here, the invariance of $h$, and the transversal invariance of $g$, we have from Lemma \ref{lemma:TransInv} that $[X,Y_0] = 0$. On the other hand, since $h$ is regular on $\operatorname{supp}(X)$ and $\dd h(Y_0) = 1$ on $V$, we have from Lemma \ref{lemma:NormalGood} that $Y_{0}$ induces a unique invariant and $h$-normalized normal class on $M$. By Theorem \ref{teo:NormalSection}, $X$ is hamiltonizable on $M$.
\end{proof}

Now, let us present some remarks on the Hamiltonization problem in the case of invariant metrics.

\begin{Lemma}\label{lemma:MetricXinv}
Let $X$ be a vector field and $g$ a transversally invariant Riemannian metric for~$X$ on an orientable manifold $M$. Then, the following assertions are equivalent:
 \begin{enumerate}\itemsep=0pt
 \item[$(a)$] The canonical volume form $\Omega_{g}$ induced by $g$ is $X$-invariant.
 \item[$(b)$] The function $g(X,X)\in C^{\infty}(M)$ is a first integral of $X$.
 \item[$(c)$] The metric $g$ is $X$-invariant.
 \end{enumerate}
\end{Lemma}
\begin{proof}
Without loss of generality, assume that $X$ is nowhere vanishing. Fix a local orthonormal basis of 1-forms $\theta_1,\dots,\theta_{m-1},\theta_{m}$ such that $\theta_1,\dots,\theta_{m-1}$ is a basis of the annihilator of $X$. Consider the matrix of functions $F=[f^{i}_{j}]$ defined by ${\rm L}_{X}\theta_{j}=\sum_{i=1}^{m}f^{i}_{j}\theta_{i}$, $j=1,\dots,m$. Note that $f^{m}_{m}=\tfrac{1}{2}\operatorname{L}_{X}\ln(g(X,X))$ and $f^{m}_{j}=0$ because $\operatorname{L}_{X}\theta_{j}$ vanishes on $X$ for $j=1,\dots,m-1$. Moreover, the transversal invariance property is equivalent to
$F = \big(\begin{smallmatrix}
 A & 0 \\
 0 & f^{m}_{m}
 \end{smallmatrix}\big)$,
where $A$ is skew-symmetric. Similarly, the $X$-invariance of $g$ is equivalent to the skew-symmetry of $F$. Finally, the canonical volume form is $\Omega_{g}=\theta_1\wedge\cdots\wedge\theta_{m}$, so ${\rm div}_{\Omega_{g}}(X)=\sum_{i=1}^{m}f^{i}_{i} = f^{m}_{m}$. Therefore, ($a$), ($b$) and ($c$) are equivalent to $f^{m}_{m}=0$.
\end{proof}

Note that, in the case of a vector field $X$ with an invariant Riemannian metric $g$, the Laplacian of a first integral is again a first integral. Indeed, let $h$ be a first integral of $X$. By Lemmas \ref{lemma:TransInv} and \ref{lemma:MetricXinv}, the gradient $\nabla h$ commutes with $X$ and the canonoical volume form $\Omega_{g}$ is $X$-invariant. Then, by Theorem \ref{teo:DiverInt}, the Laplacian $\Delta h := {\rm div}_{\Omega_{g}}(\nabla f)$ is also a first integral of $X$. In particular, there is a sequence of (not necessarily independent) first integrals $h_0,h_1,\dots$ recursively defined by $h_0:=g(X,X)$ and $h_n:=\Delta h_{n-1}$.

\begin{Corollary}
Let $X$ be a vector field, admitting an invariant Riemannian metric $g$. If~the fist integral $h_n$ is regular for some $n$, then $X$ is Hamiltonian on $M$ with respect to $\pi:=\frac{1}{g(\nabla h_{n},\nabla h_{n})}\nabla h_{n}\wedge X$ and with Hamiltonian function $h_{n}$.
\end{Corollary}

{\bf Hamiltonization of vertical vector fields.} In order to present new geometric situations in which one has positive solutions to the hamiltonization problem, we consider the case of (non necessarily regular) vector fields that are tangent to the fibers of a submersion. We begin by considering the case when the given vector field admits an invariant horizontal distribution. Also, we deal with the case when the submersion has 1-dimensional fibers. For instance, this hypothesis has been considered in \cite{FGS-05}, where it is used to give hamiltonization criteria and examples for nowhere vanishing vector fields with periodic flow.

\begin{Lemma}\label{lemma:NormalProy}
Let $p\colon M \to N$ be a submersion and $X \in \Gamma(TM)$ a vertical vector field. Let $f \in C^{\infty}(N)$ be regular on an open set $N_0$ containing $p(\operatorname{supp}(X))$, and $v \in \Gamma(TN_0)$ such that $\dd{f}|_{N_0}(v) = 1$. Then, every vector field $Y_0$ on $M_0:=p^{-1}(N_0)$ $p$-related with $v$ satisfying $[X,Y_0]\wedge X=0$ defines a $p^{*}f$-normalized and invariant normal class relative to $X$.
\end{Lemma}
\begin{proof}
Since $Y_0$ is $p$-related with $v$, we have $\dd (p^{*}f)(Y_0) = p^{*}(\dd f(v))=1$ on the open set $M_0 \supseteq \operatorname{supp}(X)$. By the fact that $[X,Y_0] \wedge X=0$ and Lemma \ref{lemma:NormalGood}, the result follows.
\end{proof}

\begin{Proposition}\label{prop:HinvHam}
Let $M \to N$ be a submersion with vertical distribution $V$ and $X \in \Gamma(TV)$. If~there exists an $X$-invariant horizontal distribution $H$, $TM = V \oplus H$, and a function $f \in C^{\infty}(N)$ regular on $p(\operatorname{supp}(X))$, then $X$ is Hamiltonian on $M$ with respect to the Poisson structure
 \begin{equation*}
 \pi = \begin{cases}
 Y_0 \wedge X & \text{on} \ \operatorname{supp}(X),\\
 0 & \text{elsewhere},
 \end{cases}
 \end{equation*}
where $Y_{0}$ is the horizontal lift of a vector field $v$ defined on the regular domain of $f$ satisfying $\dd f(v)=1$.
\end{Proposition}
\begin{proof}
Let $N_0 \subseteq N$ be the open set where $f$ is regular and $v \in \Gamma(TN_{0})$ such that $\dd{f}(v) = 1$. Define $Y_{0}:=\operatorname{hor}^H(v)$ as the horizontal lift of $v$ with respect to $H$. Since $X$ is vertical, and $Y_0$ is projectable, we have that $[X,Y_0]$ is vertical. Denote by ${\rm pr}_{H}\colon TM \to H$ the projection along the splitting $TM=V \oplus H$. Taking into account that $Y_0$ is horizontal, as well as the $X$-invariance of~$H$, we have $[X,Y_0] = [X,{\rm pr}_{H}(Y_0)] = {\rm pr}_{H}[X,Y_0] = 0$. By Lemma \ref{lemma:NormalProy} and Theorem \ref{teo:NormalSection}, we~get that $X$ is hamiltonizable on $M$.
\end{proof}

\begin{Proposition}\label{prop:HTransInv2}
Let $p\colon M \to N$ be a submersion of $1$-dimensional fibers and $X \in \Gamma(TM)$ a~ver\-ti\-cal vector field. If there exists $f \in C^{\infty}(N)$ regular on $p(\operatorname{supp}(X))$, then $X$ is Hamiltonian on $M$ with respect to the Poisson structure
 \begin{equation*}
 \pi = \begin{cases}
 Y_0 \wedge X & \text{on}  \  \operatorname{supp}(X),
 \\
 0 & \text{elsewhere},
 \end{cases}
 \end{equation*}
where $Y_{0}$ is $($any$)$ $p$-related with a vector field $v$ defined on the regular domain of $f$ satisfying $\dd f(v)=1$.
\end{Proposition}
\begin{proof}
Let $N_0 \subseteq N$ be the open set where $f$ is regular and $v \in \Gamma(TN_{0})$ such that $\dd{f}(v) = 1$. Let $M_0:=p^{-1}(N_0)$ and fix $Y_{0} \in \Gamma(TM_{0})$ $p$-related with $v$. Since $X$ is vertical, and $Y_0$ is projectable, we have that $[X,Y_0]$ is vertical. Taking into account that the $p$-fibers are 1-dimensional, we get $[X,Y_0]\wedge X=0$. By Lemma \ref{lemma:NormalProy} and Theorem \ref{teo:NormalSection}, we have that $X$ is hamiltonizable on $M$.
\end{proof}

To end this part, recall that every function on a compact manifold always has a critical point. Furthermore, the converse is also true, in the following sense \cite[Theorem 4.8]{Hirsch-61}:

\begin{Proposition}\label{prop:NonCompReg}
Every non-compact manifold admits a smooth function with no critical points.
\end{Proposition}

Then, as a direct consequence of this result, as well as of our previous propositions, we get:

\begin{Theorem}\label{teo:HorInv}
Let $M \to N$ be a submersion over a non-compact manifold $N$. Then, every vertical vector field admitting an invariant horizontal distribution is hamiltonizable on~$M$.
\end{Theorem}

\begin{Theorem}\label{teo:VertHam}
Let $M{\to}N$ be a submersion with $1$-dimensional fibers over a non-compact manifold $N$. Then, every vertical vector field is hamiltonizable on $M$.
\end{Theorem}

If $N$ is non-compact, then there exists $f\in C^{\infty}(N)$ with no critical points. Hence, Theo\-rem~\ref{teo:HorInv} follows from Proposition \ref{prop:HinvHam} and Theorem \ref{teo:VertHam} from Proposition \ref{prop:HTransInv2}.

We end this part by observing that the hamiltonization criteria of Theorems \ref{teo:IntegrableHam}, \ref{teo:RegFol} and~\ref{teo:VertHam} are equivalent. To see this, first recall that Theorem~\ref{teo:RegFol} is proven by means of Theorem~\ref{teo:IntegrableHam}. On the other hand, Theorem \ref{teo:IntegrableHam} can be seen as a particular case of Theorem \ref{teo:VertHam}. Indeed, given first integrals $h_1,\dots,h_{m-1}\in C^{\infty}(M)$ of a vector field $X$ on $M$, independent on an open dense set $U\subseteq M$, the map $p:=(h_1,\dots,h_{m-1})\colon U\subseteq M \to \R{m-1}$ is a submersion of 1-dimensional fibers such that $X$ is vertical. Moreover, the coordinate function $x_i\in C^{\infty}(\R{m-1})$ is regular and satisfies $p^*x_{i}=h_i$ and $\dd x_{i}(v_{i})=1$ for $v_i:=\frac{\partial}{\partial x_{i}}$. By Theorem \ref{teo:VertHam}, the vector field~$X$ is Hamiltonian with respect to a Poisson structure $\pi_{i}$ of rank at most 2 with Hamiltonian function~$h_{i}$. Also, it is clear that $h_{j}$ is a Casimir function of $\pi_i$ and $[\pi_{i},\pi_{j}]=0$ for $j\neq i$. Finally, we observe that Theorem~\ref{teo:VertHam} is consequence of Theorem \ref{teo:RegFol}. To see this, fix a regular function $f\in C^{\infty}(N)$ and a~vector field $v\in\Gamma(TN)$ satisfying $\dd f(v)=1$. Let $\mathcal{F}_{v}$ be the 1-dimensional foliation of $N$ by the trajectories of $v$, and $\mathcal{S}:=p^{*}\mathcal{F}_{v}$ the 2-dimensional foliation of $M$ whose leaves are the inverse images of the leaves of $\mathcal{F}_{v}$ under $p$. Note that $\mathcal{S}$ is transversal to the level sets of $h:=p^{*}f$, due to the transversality of $\mathcal{F}_{v}$ to the level sets of $f$. Moreover, since $\mathcal{S}$ contains the $p$-fibers and~$X$ is vertical, it follows that $X$ is tangent to $\mathcal{S}$. Thus, the hypothesis of Theorem~\ref{teo:RegFol} hold and~$X$ is Hamiltonian on~$M$.

The rest of this section is devoted to present some hamiltonization criteria that are motivated by Theorem~\ref{teo:VertHam}, but can be also seen as applications of Theorem~\ref{teo:TrInvM}.

{\bf The case of proper actions of 1-dimensional Lie groups.} Here we consider the case when the vector field is an infinitesimal generator of a proper action, namely, if it is complete and the action induced by its flow is proper. In this case, we can benefit from the result of Theorem~\ref{teo:VertHam}.

\begin{Theorem}\label{teo:ProperHam}
Let $X$ be a complete vector field on $M$ such that its flow induces a proper $\mathbb{R}$-action. Then, the following assertions are equivalent:
 \begin{enumerate}\itemsep=0pt
 \item[$1.$] The vector field $X$ is hamiltonizable on $M$.
 \item[$2.$] The orbit space $N := M / \mathbb{R}$ is non-compact.
 \end{enumerate}
In this case, the vector field $X$ is Hamiltonian on $M$ with respect to the Poisson structure $\pi = Y \wedge X$, where $Y \in \Gamma(TM)$ is projectable and satisfies $\dd{h}(Y)=1$, for some regular basic function $h \in C^{\infty}(M)$.
\end{Theorem}
\begin{proof}
Note that proper actions do not admit compact orbits unless the acting Lie group is compact. Thus, the proper $\mathbb{R}$-action given by the flow of $X$ is also free. So, the vector field~$X$ is nowhere vanishing and the orbit space $N$ is a smooth manifold. Therefore, we have a~one-to-one correspondence between first integrals $h \in C^{\infty}(M)$ of $X$ and functions $f \in C^{\infty}(N)$. Now, suppose that there exist $\pi$ and $h$ such that $\pi^{\sharp}\dd h = X$. Then, $h$ is a first integral of $X$ with no critical points and the corresponding $f \in C^{\infty}(N)$ has no critical points. Thus, $N$ cannot be compact. Conversely, if $N$ is non-compact, then the result follows from Theorem \ref{teo:VertHam}.
\end{proof}

Now, we consider the class of vector fields with periodic flow such that their orbits are contained in the fibers of an $\mathbb{S}^1$-bundle. This includes vector fields that differ by a first integral scalar factor from the infinitesimal generator of an $\mathbb{S}^{1}$-action. In particular:

\begin{Lemma}\label{lemma:Periodic}
 The trajectories of a nowhere vanishing periodic vector field are the fibers of an $\mathbb{S}^1$-bundle.
\end{Lemma}
\begin{proof}
Let $X$ be a nowhere vanishing vector field on $M$ with periodic flow, $\varpi\in C^{\infty}(M)$ the period of $X$ and $\Upsilon:=1/\varpi X$. By construction, $\Upsilon$ has period $2\pi$, so it is the infinitesimal generator of a free $\mathbb{S}^1$-action on $M$ whose orbits agree with the trajectories of $X$.
\end{proof}

With these ingredients we formulate the following:

\begin{Theorem}\label{teo:Periodic}
Let $X$ be a nowhere vanishing vector field with periodic flow on $M$. Then, $X$ is hamiltonizable on $M$ if and only if $M$ is non-compact. In this case, $M$ admits an $\mathbb{S}^1$-bundle structure and the vector field $X$ is Hamiltonian on $M$ with respect to the Poisson structure $\pi = Y \wedge X$, where $Y$ is any projectable vector field on $M$ satisfying $\dd{h}(Y)=1$, for some regular basic function $h \in C^{\infty}(M)$.
\end{Theorem}
\begin{proof}
If $M$ is compact, then any function $h \in C^{\infty}(M)$ has a critical point. Since $X$ is nowhere vanishing, the equation $\pi^{\sharp}\dd h = X$ cannot globally hold on $M$. Conversely, if $M$ is non-compact, then by Lemma \ref{lemma:Periodic} we get
that $X$ is a vertical vector field for some $\mathbb{S}^1$-bundle $M \to N$. Since~$M$ is non-compact, so is not $N$. By Theorem \ref{teo:VertHam}, $X$ is hamiltonizable.
\end{proof}

The result of Theorem \ref{teo:Periodic}, is related to the hamiltonization criteria and examples given in \cite{FGS-05}. Indeed, Theorem~\ref{teo:Periodic} characterizes the solvability of the hamiltonization problem for non-vanishing vector fields with periodic flow, improving \cite[Theorem~1]{FGS-05}. More precisely, we have shown that the \emph{fibrating-periodic flow hypothesis} is always satisfied, due to Lem\-ma~\ref{lemma:Periodic}. Furthermore, by means of Proposition~\ref{prop:NonCompReg}, the existence of the regular first integral is translated into the non-compacity of the manifold.

Finally, we present the following hamiltonization criteria for non-necessarily nowhere vani\-shing vector fields with periodic flow.

\begin{Theorem}\label{teo:Periodic2}
Every vector field $X$ on $M$ with periodic flow admitting a first integral $h$ that is regular on $\operatorname{supp}(X)$ is Hamiltonian on $M$ with respect to the Poisson structure $\pi = Y \wedge X$, where $Y = \big\langle \widetilde{Y} \big\rangle$ is the averaging of a vector field $\widetilde{Y}$ satisfying $\dd h\big(\widetilde{Y}\big)X = X$ with respect to the $\mathbb{S}^1$-action induced by the flow of~$X$.
\end{Theorem}
\begin{proof}
By Lemma \ref{lemma:NormalGood}, there exists $\widetilde{Y} \in \Gamma(TM)$ such that $\dd h\big(\widetilde{Y}\big)X = X$. Now, consider the average $Y := \big\langle \widetilde{Y} \big\rangle$ of $\widetilde{Y}$ under the $\mathbb{S}^{1}$-action induced by the flow of $X$. Then, we have $[X,Y]=fX$, for some $f \in C^{\infty}(M)$. Moreover, by the $\mathbb{S}^{1}$-invariance of $h$ and $X$, and the properties of the averaging operator, we get that $\dd h(Y)X = X$. In other words, the normal class induced by $Y$ is invariant and $h$-normalized. By Theorem~\ref{teo:NormalSection}, the vector field $X$ is hamiltonizable on~$M$.\looseness=-1
\end{proof}

Observe that the results of this part can be also obtained by means of Theorem~\ref{teo:TrInvM}. Indeed, the Poisson structures in Theorems \ref{teo:ProperHam}, \ref{teo:Periodic} and \ref{teo:Periodic2} can alternatively be constructed as
 \begin{equation*}
 \pi=\frac{\nabla h}{g(\nabla h,\nabla h)}\wedge X,
 \end{equation*}
where $g$ is any $X$-invariant Riemmanian metric and $h \in C^{\infty}(M)$ is a regular first integral of~$X$. The existence of $g$ follows from the properness of the corresponding Lie group action \cite[Proposition 2.5.2]{Duistermaat-1999}, and the existence of $h$ in Theorems \ref{teo:ProperHam} and \ref{teo:Periodic} follows from Proposition~\ref{prop:NonCompReg} and the non-compactness of the orbit space. Finally, the fact that $X$ is Hamiltonian with respect to~$\pi$ follows from Theorem \ref{teo:TrInvM}.

We now illustrate our hamiltonization criteria for periodic vector fields (compare with Example \ref{ex:TorusIntegrable}).

\begin{Example}
Consider the torus $\mathbb{T}^2$ with natural coordinates $(\varphi_1,\varphi_2)$, $\varphi_i \in \mathbb{R}/2\pi\mathbb{Z}$.
Recall that, for coprime integers $m$ and $n$, the vector field $\Upsilon = m \tfrac{\partial}{\partial\varphi_1} + n\tfrac{\partial}{\partial\varphi_2}$ is the infinitesimal generator of an $\mathbb{S}^1$-action on $\mathbb{T}^2$. Now, fix $2\pi$-periodic functions $F,\omega\in C^{\infty}(\mathbb{R})$, and define
 \begin{equation*}
 X := \omega(n\varphi_1-m\varphi_2) \Upsilon,
 \end{equation*}
and $h(\varphi_{1}, \varphi_{2}):=F(n\varphi_{1} - m\varphi_{2})$. Then, the function $h$ is a first integral of $X$.
Moreover, for integers $r$ and $s$ satisfying $nr-ms=1$, the vector field
 \begin{equation*}
 Y = \frac{r}{F'(n\varphi_{1} - m\varphi_{2})}\frac{\partial}{\partial\varphi_1} + \frac{s}{F'(n\varphi_{1} - m\varphi_{2})}\frac{\partial}{\partial\varphi_2}
 \end{equation*}
defines an invariant and $h$-normalized normal class relative to $X$ on its domain.
Therefore, $X$ is a Hamiltonian vector field with Hamiltonian function $h$ and with respect to the Poisson structure
 \begin{equation*}
 \pi := Y \wedge X = \frac{\omega(n\varphi_1-m\varphi_2)}{F'(n\varphi_{1} - m\varphi_{2})}\, \frac{\partial}{\partial \varphi_1}\wedge \frac{\partial}{\partial \varphi_2}.
 \end{equation*}
\end{Example}

{\bf Hamiltonization of infinitesimal generators.} Here, we adapt the result of Theorem \ref{teo:TrInvM} to obtain a hamiltonization criteria for infinitesimal generators of proper actions.

\begin{Proposition}\label{prop:InfHam}
Let $G$ be a Lie group with Lie algebra $\mathfrak{g}$, acting properly on $M$. For each $\xi \in \mathfrak{g}$, let $\xi_{M}$ be the infinitesimal generator of the $G$-action.
Suppose that there exists a $G$-invariant function $h \in C^{\infty}(M)$ with no critical points.
Then, for every $G$-invariant Riemannian metric $g$, we have the linear map
 \begin{equation*}
 \xi \mapsto \pi_{\xi}:=\frac{\nabla h}{g(\nabla h,\nabla h)}\wedge\xi_{M}
 \end{equation*}
from $\mathfrak{g}$ to a vector space of Poisson structures of rank at most two on $M$ such that $\pi_{\xi}^{\sharp}\dd h = \xi_{M}$.
\end{Proposition}
\begin{proof}
Since $G$ acts properly on $M$, there exists an invariant Riemannian metric $g$ on $M$ \cite[Proposition 2.5.2]{Duistermaat-1999}. Now, consider the vector field $Y := \tfrac{\nabla h}{g(\nabla h,\nabla h)}$, where $\nabla h$ denotes the gradient vector field of $h$. Clearly, we have $\dd h(Y) = 1$. On the other hand, since $h$ and $g$ are invariant, we get that $Y$ is invariant. In other words, for each $\xi \in \mathfrak{g}$, we have $[\xi_{M},Y] = 0$ and $\dd h(\xi_{M}) = 0$. By setting $\pi_{\xi} := Y \wedge \xi_{M}$, the result follows.
\end{proof}

We also have the following hamiltonization criteria for actions that are proper and free.

\begin{Proposition}\label{prop:InfHam2}
Let $G$ be a Lie group acting freely and properly on $M$. If the orbit space $M/G$ is non-compact, then there exists a linear map $\xi \mapsto \pi_{\xi}:=\frac{\nabla h}{g(\nabla h,\nabla h)}\wedge\xi_{M}$ from the Lie algebra $\mathfrak{g}$ of $G$ to a vector space of Poisson structures of rank at most two on $M$ such that $\pi_{\xi}^{\sharp}\dd h = \xi_{M}$.
\end{Proposition}
\begin{proof}
By Proposition \ref{prop:InfHam}, it suffices to show that there exists a regular $G$-invariant function $h$ on $M$. Since $N$ is non-compact, we have from Proposition \ref{prop:NonCompReg} that there exists $f \in C^{\infty}(M/G)$ with no critical points. If $p\colon M\to M/G$ is the canonical projection to the orbit space, then the pull-back function $h:=p^{*}f$ is $G$-invariant, and also regular due to the regularity of $f$.
\end{proof}

\begin{Remark}
Note that Proposition \ref{prop:InfHam2} also follows from Theorem \ref{teo:HorInv}. Indeed, an invariant horizontal distribution always exists because of the $G$-action is free and proper.
\end{Remark}

We now illustrate the previous ideas on the vector fields of Example \ref{ex:Torus1}.

\begin{Example}\label{ex:Torus2}
Consider the vector fields on $\R{3}$
 \begin{equation*}
 X_{1} = 2xz\frac{\partial}{\partial x} + 2yz\frac{\partial}{\partial y} + \big(1-x^2-y^2+z^2\big)\frac{\partial}{\partial z} \qquad \text{and} \qquad
 X_{2} = y\frac{\partial}{\partial x} - x\frac{\partial}{\partial y}.
 \end{equation*}
Note that $X_1$ and $X_2$ commute and have periodic flow. Therefore, the flows of $X_1$ and $X_2$ induce a $\mathbb{T}^2$-action on $\R{3}$, which is free on $M = \R{3} \setminus \big(\{z-\text{axis}\}\cup\big\{x^2+y^2=1, z = 0\big\}\big)$. Indeed, $X_1$ and $X_2$ have as common first integral $h = \frac{x^2+y^2}{(x^2+y^2+z^2+1)^2}$, whose level sets on $M$ are the orbits of the $\mathbb{T}^2$-action. Now, let $Y_{0}$ any vector field on $M$ such that $\dd h(Y_{0}) = 1$. Since $h$ is $\mathbb{T}^2$-invariant, the averaged vector field $Y := \langle Y_{0}\rangle^{\mathbb{T}^2}$ is invariant and also satisfies $\dd h(Y) = 1$. Note that for each $\lambda \in \mathbb{R}$, the Poisson structure
 \begin{equation*}
 \pi_{\lambda} := Y\wedge X_{\lambda}, \qquad \text{where}\quad X_{\lambda} := X_1 + \lambda X_2,
 \end{equation*}
is such that $\pi_{\lambda}^{\sharp}\dd h = X_{\lambda}$. In other words, the vector field $X_{\lambda}$ is hamiltonizable by Proposition~\ref{prop:InfHam}, regardless of the value of $\lambda$. Moreover, recall from Example \ref{ex:Torus1} that if $\lambda$ is irrational, then the orbits of $X_{\lambda}$ are dense on the level sets of $h$, which are compact. Thus, its flow does not define a proper action, and therefore the vector field $X_{\lambda}$ cannot be hamiltonized by the result of Theorem \ref{teo:ProperHam}.
\end{Example}

 \section{Hamiltonization of torus actions}

In this section we address the question on hamiltonization of Lie group actions for the particular case in which the acting Lie group is a torus. We adapt some of the ideas developed in Section~\ref{sec:HamGood} to provide a Poisson structure so that a given torus action is Hamiltonian with momentum map.

Let $X_{1}, \dots, X_{k} \in \Gamma(TM)$ be infinitesimal generators of an action of $\mathbb{T}^{k}$ on $M$. Suppose also that we are given $\mathbb{T}^{k}$-invariant functions $h_{1},\dots,h_{k} \in C^{\infty}(M)$ satisfying the following condition: there exist $Y_1,\dots,Y_{k} \in \Gamma(TM)$ such that
 \begin{equation}\label{eq:Goodness}
 \sum_{j=1}^{k}\dd h_{i}(Y_{j})X_{j} = X_{i},\qquad \text{for every}\quad i \in \{1,\dots,k\}.
 \end{equation}

\begin{Remark}
In the case $k=1$, condition \eqref{eq:Goodness} means that $h_{1}$ is a first integral of $X_{1}$ that is regular on $\operatorname{supp}(X_{1})$, due to Lemma \ref{lemma:NormalGood}.
\end{Remark}

\begin{Remark}
By the compactness of $\mathbb{T}^k$, the action admits an invariant Riema\-niann metric~$g$. Then, from Lemma \ref{lemma:MetricXinv}, in the case when $X_1, \dots, X_{k}$ are nowhere vanishing, we~have the following family of first integrals $h_{n;i,j}:=\Delta^{n-1}g(X_i,X_j)$, $i,j=1,\dots,k$ and $n\in\mathbb{N}$, where $\Delta$ is the Laplacian of $g$, from which we may be able to choose~$k$ of them satisfying \eqref{eq:Goodness}.
\end{Remark}

Let also $\R{k} = {\rm Lie}\big(\mathbb{T}^{k}\big)$ be the abelian Lie algebra of $\mathbb{T}^{k}$. For each $\xi = (\xi_1,\dots,\xi_{k}) \in \R{k}$, consider the infinitesimal generator $\xi_{M} := \sum_{i=1}^{k}\xi_{i}X_{i}$ and $h_{\xi} \in C^{\infty}(M)$ given by $h_{\xi} := \sum_{i=1}^{k}\xi_{i}h_{i}$.

Now, define $\pi \in \Gamma\big({\wedge}^{2}TM\big)$ by
 \begin{equation}\label{eq:PoissTorus}
 \pi := \sum_{j=1}^{k}Y_{j}\wedge X_{j}.
 \end{equation}

\begin{Lemma}\label{lemma:PoissToric0}
For each $\xi \in \R{k}$, the identity $\pi^{\sharp}\dd h_{\xi} = \xi_{M}$ holds. In particular, in the case when~$\pi$ is Poisson, the $\mathbb{T}^{k}$-action on $M$ is Hamiltonian with momentum map $\mathbb{J}\colon M \to \big(\R{k}\big)^{*}$ given by $\mathbb{J}(x)(\xi) := h_{\xi}(x)$.
\end{Lemma}
\begin{proof}
Since $h_{i}$ is a first integral of $X_{i}$, we have from the condition \eqref{eq:Goodness} that $\pi^{\sharp} \dd h_{i} = \sum_{j=1}^{k} \dd h_{i} (Y_{j}) X_{j} = X_{i}$. Therefore, $\pi^{\sharp}\dd h_{\xi} = \sum_{i=1}^{k}\xi_{i}\pi^{\sharp}\dd h_{i} = \sum_{i=1}^{k}\xi_{i}X_{i} = \xi_{M}$.
\end{proof}

Now, let us study conditions under which $\pi$ is a Poisson structure. First, observe that, since the action is abelian, the infinitesimal generators $X_{i}$ commute, $[X_{i},X_{j}] = 0$. Moreover, by averaging with respect to the $\mathbb{T}^{k}$-action, and taking into account the invariance of each~$h_i$ and~$X_i$, one can assume that the vector fields $Y_{i}$ given in the condition \eqref{eq:Goodness} are also invariant, $[X_{i},Y_{j}] = 0$. Therefore,
 \begin{equation*}
 [\pi,\pi] = \sum_{1 \leq i,j \leq k}[Y_i,Y_{j}]\wedge X_{i}\wedge X_{j}.
 \end{equation*}

\begin{Lemma}\label{lemma:PoissToric1}
Let $U_{ij} \subseteq M$ be the open set on which $X_{i} \wedge X_{j} \neq 0$. If, for every $i,j \in \{1,\dots,k\}$, one has $[Y_{i},Y_{j}]|_{U_{ij}} = a_{ij}X_{i} + b_{ij}X_{j}$, for some $a_{ij},b_{ij} \in C^{\infty}(U_{ij})$, then $\pi$ is a Poisson structure on $M$.
\end{Lemma}

{\bf Actions of 2-dimensional tori.} Recall that the Poisson property for $\pi$ in \eqref{eq:PoissTorus} holds in particular for $k=1$: this is the content of Theorem \ref{teo:Periodic2}. Now, we focus on the case of an action of a 2-dimensional torus on $M$, with infinitesimal generators $X_1$, $X_2$ and invariant functi\-ons~$h_1$,~$h_2$ satisfying the condition
 \begin{equation}\label{eq:Goodness2}
 \dd h_{i}(Y_{1})X_{1} + \dd h_{i}(Y_{2})X_{2} = X_{i}, \qquad i = 1, 2,
 \end{equation}
for some $\mathbb{T}^2$-invariant $Y_1,Y_2 \in \Gamma(TM)$. We then have that the bivector field
 \begin{equation}\label{eq:PoissTorus2}
 \pi = Y_1\wedge X_1 + Y_2\wedge X_2
 \end{equation}
satisfies $[\pi,\pi] = [Y_1,Y_2] \wedge X_1 \wedge X_2$.
Now, let $U$ be the open set in which $X_{1}$ and $X_{2}$ are independent.

\begin{Lemma}\label{lemma:PoissToric2}
For $i = 1,2$, we have $\dd h_{i}[Y_{1},Y_{2}]|_{U} = 0$. In particular, $\ii_{\dd h_{i}}[\pi,\pi] = 0$.
\end{Lemma}
\begin{proof}
By \eqref{eq:Goodness2}, ${\rm L}_{Y_{i}}h_{j} = \dd h_{j}(Y_{i}) = \delta_{ij}$. Thus, $\dd h_{i}[Y_{1},Y_{2}] = {\rm L}_{Y_{1}}{\rm L}_{Y_{2}}h_{i} - {\rm L}_{Y_{2}}{\rm L}_{Y_{1}}h_{i} = 0$ on $U$.
\end{proof}

\begin{Proposition}\label{prop:2Toric}
Suppose that $\operatorname{dim}{M} = 4$. Then, the bivector field $\pi$ in \eqref{eq:PoissTorus2} is Poisson on~$M$ and has rank four on $U$.
\end{Proposition}
\begin{proof}
Condition \eqref{eq:Goodness2} implies that $h_{1}$ and $h_{2}$ are functionally independent on $U$, so their common level sets induce a 2-dimensional foliation $\mathcal{L}$ of $U$. Moreover, by the invariance of $h_{1}$ and $h_{2}$, the independent vector fields $X_{1}$ and $X_{2}$ span $T\mathcal{L}$. On the other hand, by Lemma \ref{lemma:PoissToric2}, the vector field $[Y_1,Y_2]$ is tangent to $\mathcal{L}$ on $U$. In particular, the hypothesis in Lemma \ref{lemma:PoissToric1} is satisfied, proving that $\pi$ is Poisson on $M$. Finally, from the condition \eqref{eq:Goodness2}, the vector fields $Y_1$ and $Y_2$ are independent and normal to $\mathcal{L}$ on $U$. Therefore, the vector fields $Y_1$, $Y_2$, $X_1$, $X_2$ are independent and $\pi$ has maximal rank on $U$.
\end{proof}

In general, the orbits of an action are contained in the level sets of invariant functions. In~Pro\-position \ref{prop:2Toric}, the dimension hypothesis allows to describe the orbits as level sets of invariant functions. This motives the consideration of the following setting, in which the dimension hypothesis can be relaxed: the orbits and the level sets agree along the leaves of a foliation of the manifold, with an additional transversality condition.

Denote by $\mathcal{L}$ the foliation of $M$ given by the common level sets of $h_1$ and $h_2$, and by $\mathcal{O}$ the foliation of $M$ by the orbits of the $\mathbb{T}^2$-action.

\begin{Definition}\label{def:compat}
A foliation $\mathcal{S}$ of $M$ is said to be \emph{compatible} with $\mathcal{L}$ and $\mathcal{O}$ if the following conditions are satisfied:
 \begin{enumerate}\itemsep=0pt
 \item The foliations $\mathcal{S}$ and $\mathcal{L}$ are transversal, $T\mathcal{S} + T\mathcal{L} = TM$.
 \item The leaves of the foliation $\mathcal{S}$ are invariant, $T\mathcal{O} \subset T\mathcal{S}$.
 \item On the open set $U$ in which $X_1 \wedge X_2 \neq 0$, one has $T\mathcal{O} = T\mathcal{S} \cap T\mathcal{L}$.
 \end{enumerate}
\end{Definition}

Observe that, in the context of Proposition \ref{prop:2Toric}, such a compatible foliation exists: it consists of the connected components of $M$. More generally, we have:

\begin{Proposition}\label{prop:PoissTorus}
The bivector field $\pi$ in \eqref{eq:PoissTorus2} is Poisson on $M$ if and only if there exists a~foli\-a\-tion $\mathcal{S}$ on $M$ that is compatible with $\mathcal{L}$ and $\mathcal{O}$ on $U$ and $Y_1$, $Y_2$ are tangent to $\mathcal{S}|_{U}$.
\end{Proposition}
\begin{proof}
By condition \eqref{eq:Goodness2}, the invariant functions $h_1$ and $h_2$ are independent on $U$, as well as $Y_1$ and $Y_2$. In particular, the foliation $\mathcal{L}$ is regular on $U$. Now, suppose that there exists a~foliation~$\mathcal{S}$ on $M$ such that is compatible with $\mathcal{L}$ and $\mathcal{O}$ on $U$ and $Y_1,Y_2 \in \Gamma(T\mathcal{S}|_{U})$. Then, $[Y_1,Y_2]|_{U} \in \Gamma(T\mathcal{S})$. On the other hand, by Lemma~\ref{lemma:PoissToric2}, we have $\dd h_{i}[Y_1,Y_2]|_{U} = 0$, so $[Y_1,Y_2]|_{U} \in \Gamma(T\mathcal{L})$. By the compatibility of $\mathcal{S}$, this implies that $[Y_1,Y_2]|_{U} \in \Gamma(T\mathcal{O})$. Finally, since~$X_1$ and~$X_2$ generate~$T\mathcal{O}$, we have $[Y_1,Y_2] \wedge X_1 \wedge X_2|_{U} = 0$. From here, and the fact that $X_1 \wedge X_2 = 0$ outside of $U$, we conclude that $\pi$ is Poisson on $M$. Conversely, suppose that~$\pi$ is Poisson. Let us show that the symplectic foliation $\mathcal{S}$ of $\pi$ is compatible with $\mathcal{L}$ and $\mathcal{O}$ on $U$, and satisfies $Y_1,Y_2 \in \Gamma(T\mathcal{S}|_{U})$. By \eqref{eq:Goodness2}, we have $\dd h_{j}(X_i) = 0$ and $\dd h_{j}(Y_i) = \delta_{ji}$ on~$U$, so the vector fields $Y_1$, $Y_2$, $X_1$, $X_2$ are independent on $U$. Thus, there exist $\alpha_1,\alpha_2 \in \Gamma(T^{*}U)$ such that $\alpha_{i}(X_{j}) = -\delta_{ij}$ and $\alpha_{i}(Y_{j}) = 0$. Then, $\pi^{\sharp}\alpha_{i} = Y_{i}$, so $Y_1$ and $Y_2$ are tangent to $\mathcal{S}$ on $U$. From here, we also have that $\mathcal{S}$ and $\mathcal{L}$ are transversal on $U$:
 \begin{equation*}
 T\mathcal{S}|_{U}+T\mathcal{L}|_{U} \supseteq {\rm span}\{Y_1|_{U},Y_2|_{U}\} + T\mathcal{L}|_{U} = TM|_{U}
 \end{equation*}
Finally, if $v \in (T\mathcal{S} \cap T\mathcal{L})|_U$, then $v = a_1 X_1+a_2X_2 + b_1Y_2+b_2Y_2$, with $b_{i} = \dd h_{i}(v) = 0$. So, $v \in T\mathcal{O}$, showing that $T\mathcal{S} \cap T\mathcal{L}=T\mathcal{O}$.
\end{proof}

We now give sufficient conditions under which the vector fields $Y_{1}$ and $Y_{2}$ in \eqref{eq:Goodness2} can be chosen as sections of a given compatible foliation.

\begin{Lemma}\label{lemma:Compat}
Let $V$ be a $\mathbb{T}^{2}$-saturated open set satisfying $\operatorname{supp}(X_1) \cup \operatorname{supp}(X_2) \subseteq V$ and $\dd h_1 \wedge \dd h_2 \neq 0$ on $V$. If $\mathcal{S}$ is a compatible regular foliation on $V$, then there exist $\mathbb{T}^{2}$-invariant vector fields $Y_1, Y_{2} \in \Gamma(TM)$ tangent to $\mathcal{S}$ on $V$ and satisfying \eqref{eq:Goodness2}.
\end{Lemma}
\begin{proof}
Since $\dd h_{1}$ and $\dd h_{2}$ are invariant and independent in the saturated open set $V$, there exist invariant $\widetilde{Y}_{1},\widetilde{Y}_{2} \in \Gamma(TV)$ satisfying $\dd h_{i}\big(\widetilde{Y}_{j}\big) = \delta_{ij}$. Since $\mathcal{S}$ and $\mathcal{L}$ are regular and transversal on~$V$, there exist $Y'_1,Y'_2 \in \Gamma(T\mathcal{S})$ and $Z'_1,Z'_2 \in \Gamma(T\mathcal{L})$ such that $\widetilde{Y}_{j} = Y'_{j} + Z'_{j}$. Furthermore, since $T\mathcal{O} \subseteq T\mathcal{S} \cap T\mathcal{L}$, we have that the action restricts to the leaves of these foliations. So, by averaging and taking into account the invariance of $\widetilde{Y}_{i}$, we may assume that $Y'_{j}$, $Z'_{j}$ are invariant. Since $\dd h_{i}(Z'_{j})=0$, we have $\dd h_{i}(Y'_{j}) = \delta_{ij}$. On the other hand, observe that the closed sets $C := \operatorname{supp}(X_1) \cup \operatorname{supp}(X_2)$ and $M \setminus V$ are disjoint. Thus, there exist disjoint open sets $V'$ and $W'$ such that $C \subseteq W'$ and $M \setminus V \subseteq V'$. Let $\mu \in C^{\infty}(M)$ be such that $\mu|_{C} = 1$ and $\mu|_{\overline{V'}} = 0$. Since $X_1$, $X_2$ are trivial outside $C$ and $\mu$ is constant on $C$, we have that $\mu$ is invariant. Then, there exist smooth, global and invariant vector fields $Y_1,Y_2 \in \Gamma(TM)$, well-defined by $Y_{j} := \mu Y'_{j}$ on $V$ and $Y_{j} = 0$ on $V'$. It is left to show that $Y_1$, $Y_2$ satisfy the condition \eqref{eq:Goodness2}. First, observe that the condition is automatically satisfied on $M \setminus C$. On the other hand, since $Y_{j}$ agrees with~$Y'_{j}$ on $C \subseteq V$, we get that $\dd h_{i}(Y_{j}) = \delta_{ij}$ on $C$. Therefore, the condition is satisfied.
\end{proof}

{\sloppy\begin{Theorem}\label{teo:2Torus}
Let $X_1$ and $X_{2}$ be infinitesimal generators of a $\mathbb{T}^{2}$-action on a~manifold~$M$, and $h_{1},h_{2} \in C^{\infty}(M)$ $\mathbb{T}^{2}$-invariant functions independent on a saturated open neighborhood $V$ of $\operatorname{supp}(X_1) \cup \operatorname{supp}(X_2)$. Suppose that there exists a compatible regular foliation $\mathcal{S}$ on $V$, in the sense of Definition $\ref{def:compat}$. Then, there exists a Poisson structure $\pi$ on $M$ such that the $\mathbb{T}^{2}$-action is Hamiltonian on $(M,\pi)$ with momentum map $\mathbb{J}\colon M \to (\R{2})^{*}$ given by $\mathbb{J}(x)(\xi) := \xi_1 h_1(x) + \allowbreak\xi_2h_2(x)$, $\xi \in \R{2}$.
\end{Theorem}}

\begin{proof}
Let $U$ be the open set in which $X_1 \wedge X_2\neq 0$. By Lemma \ref{lemma:Compat}, there exist invariant vector fields $Y_1$, $Y_2$ on $M$ tangent to $\mathcal{S}$ on $V$ satisfying \eqref{eq:Goodness2}. Since $\mathcal{S}$ is compatible on $U$, the bivector field $\pi$ in \eqref{eq:PoissTorus2} is Poisson, by Proposition \ref{prop:PoissTorus}. Finally, by Lemma \ref{lemma:PoissToric0}, we get that $\mathbb{J}$ is a momentum map.
\end{proof}

\begin{Remark}
We believe that Lemma \ref{lemma:Compat}, and hence Theorem \ref{teo:2Torus}, are still true in the more general case when the invariant functions $h_1$, $h_2$ satisfy the condition \eqref{eq:Goodness2}, not only when $\operatorname{supp}(X_1) \cup \operatorname{supp}(X_2) \subseteq V$. However, establishing such fact requires a more profound analysis of the condition \eqref{eq:Goodness2}, perhaps by providing an intrinsic formulation of it.
\end{Remark}

\subsection*{Acknowledgements}

We are very grateful to the anonymous referees for the observations and suggested improvements on various aspects of this work.
This research was partially supported by the Mexican National Council of Science and Technology (CONACYT) under the grant CB2015 no.~258302 and the University of Sonora (UNISON) under the project no.~USO315007338.
J.C.R.P.\ thanks CONACyT for a postdoctoral fellowship held during the production of this work.
E.V.B.\ was supported by FAPERJ grants E-26/202.411/2019 and E-26/202.412/2019.

\pdfbookmark[1]{References}{ref}
\LastPageEnding

\end{document}